\documentclass[a4paper,reqno]{amsart}
\usepackage{pdfsync}
\usepackage{amsmath,amssymb}
\usepackage{enumerate}
\usepackage{xspace}
\usepackage{boxedminipage}
\usepackage{algorithmic}
\usepackage{listings}
\usepackage{tabularx}
\usepackage{subfigure}

\newenvironment{Proof}[1][{. }]
	       {\noindent{\bf \Proofname\ #1}}
	       {{\raggedright{{ }\hfill\qed}}} 
	 
\newtheoremstyle%
    {plain}
    {}
    {}
    {\mdseries\slshape}
    {}
    {\bfseries}
    {.}
    {1.0ex}
    {}
    
\newtheoremstyle
    {note}
    {}
    {}
    {}
    {}
    {\bfseries}
    {.}
    {1.0ex}
    {}
    
\swapnumbers
\theoremstyle{plain}
\newtheorem{The}[subsection]{Theorem}

\newtheorem{Pro*}{Proposition}
\newtheorem{Lem}[subsection]{Lemma}

\theoremstyle{note}
\newtheorem{Example}[subsection]{Example}
\newtheorem{Rem}[subsection]{Remark}
\newtheorem{Obs}[subsection]{Remark}
\newtheorem{Defn}[subsection]{Definition}

\usepackage{ifthen}
\usepackage{amsmath,amssymb}    
\usepackage{stmaryrd} 
\usepackage{eufrak}     
\usepackage{mathrsfs} 
\usepackage{bbold}    
\provideboolean{usemathrsfs}
\setboolean{usemathrsfs}{true}
\usepackage{enumerate}
\usepackage{xspace}
\usepackage{boxedminipage}
\usepackage[latin1]{inputenc}
\usepackage{verbatim}
\provideboolean{usebeamer}
\provideboolean{isthesis}
\provideboolean{isamsltex}
\provideboolean{issiamltex}

\usepackage{hyperref}
\newcommand{\linkedurl}[1]{{\color{blue}\href{#1}{#1}}}
\newcommand{\linkedemail}[1]{{\makeatletter\color{blue}\href{mailto:#1}{#1}\makeatother}}
\newcommand{\Ignore}[1]{}
\newcommand{\freeze}[1]{}
\newcommand{\crossout}[1]{{\textcolor{red}{#1}}}
\newcommand{\highlight}[1]{{\textcolor{blue}{#1}}}

\provideboolean{shownotes}
\setboolean{shownotes}{true}
\newcounter{margnote}[page]
\newcommand{\margnotemark}{\highlight{\upshape\texttt{>\arabic{margnote}<}}}
\newcommand{\margnote}[2][]{\ifthenelse{\boolean{shownotes}}%
{\stepcounter{margnote}%
\margnotemark%
\marginpar{\texttt{\raggedright\tiny\margnotemark{#1}: #2}}}
{}}
\provideboolean{showtodo}

\newcommand{\todo}[1]
{\ifthenelse{\boolean{showtodo}}{\margnote[To do.]{#1}}{}}
{\ifthenelse{\boolean{showtodo}}{\end{boxedminipage}}{}}

\provideboolean{showcomments}
\newcommand{\margincomment}[1]{
\ifthenelse{\boolean{showcomments}}{\marginpar{\tiny #1}}{}
}
\provideboolean{showchanges}
\setboolean{showchanges}{false}
\newcommand{\changes}[1]{
  \ifthenelse{\boolean{showchanges}}
	     {{\highlight{#1}}}
	     {#1}
}
\newcommand{\changefromto}[2]{
  \ifthenelse{\boolean{showchanges}}
  {{\crossout{#1}${\color{magenta}\mapsto}$}{\highlight{#2}}}
  {#2}
}

\usepackage{mathrsfs}
\provideboolean{usemathrsfs}
\setboolean{usemathrsfs}{true}
\newcommand{\mathscript}
	   {\mathscr}

 \newcommand{\cL}{\ensuremath{\mathscript L}\xspace}

 \newcommand{\cT}{\ensuremath{\mathscript T}\xspace}

 \newcommand{\cX}{\ensuremath{\mathscript X}\xspace}
 \newcommand{\cY}{\ensuremath{\mathscript Y}\xspace}

 \newcommand{\rN}{\ensuremath{\mathbb N}\xspace}
 
 \newcommand{\rP}{\ensuremath{\mathbb P}\xspace}
 
 \newcommand{\rR}{\ensuremath{\mathbb R}\xspace}

 \newcommand{\rX}{\ensuremath{\mathbb X}\xspace}
 \newcommand{\rY}{\ensuremath{\mathbb Y}\xspace}

 \newcommand{\naturals}{\rN\xspace}

 \newcommand{\reals}{\rR}
 \newcommand{\R}[1]{\reals^{#1}}

 \newcommand{\closure}[1]{\overline{#1}}

 \newcommand{\W}{\ensuremath{\varOmega}\xspace}
 
 \renewcommand{\a}{\ensuremath{\alpha}\xspace}
 \renewcommand{\b}{\ensuremath{\beta}\xspace}

 \newcommand{\w}{\ensuremath{\omega}\xspace}

 \newcommand{\qp}[1]{\ensuremath{\left({#1}\right)}}

 \newcommand{\powqp}[2]{\ensuremath{\qp{#2}^{\kern -.2em\lower .7ex\hbox{\scriptsize $#1$}}\kern-.3em}}
 
 \newcommand{\norm}[1]{\ensuremath{\left|#1\right|}}
 
 \newcommand{\Norm}[1]{\ensuremath{\left\|#1\right\|}}
 
 \newcommand{\ltwop}[2]{\ensuremath{\left\langle#1,#2\right\rangle}}

 \newcommand{\duality}[2]{\ensuremath{\left\langle #1\,\vert\,#2\right\rangle}}
 
 \newcommand{\ensemble}[2]{\ensuremath{\left\{ #1:\;#2 \right\}}}

 \newcommand{\rangefromto}[3]{\ensuremath{#1\in\fromto{#2}{#3}}}

 \providecommand{\D}{\ensuremath{\mathrm{D}}}

 \newcommand{\registered}%
	    {\ensuremath{{}^{\bigcirc\!\;\!\!\!\!\!\!\!\;\text{\sc r}}}}

\newcommand{\AND}{\ensuremath{\text{ and }}}

\renewcommand{\div}{\operatorname{div}}

\newcommand{\Oh} {\operatorname{O}}                   
\newcommand{\pd}[2]{\ensuremath{\partial_{#1}{#2}}\xspace} 

\newcommand{\transpose}{{\boldsymbol\intercal}}   
\newcommand{\Transpose}[1]{\ensuremath{{#1}^{\transpose}}}

\newcommand{\Hess}{\ensuremath{\D^2}}
\newcommand{\intersected}{\ensuremath{\cap}}
\newcommand{\meet}{\intersected}

\newcommand{\union}[1]{\ensuremath{\bigcup}_{#1}}

\renewcommand{\vec}[1]{\ensuremath{\boldsymbol{#1}}}
\newcommand{\geovec}[1]{\ensuremath{\vec{#1}}}

\newcommand{\mat}[1]{\vec{#1}}
\newcommand{\num}[1]{\ensuremath{\mathsf{#1}}}
\newcommand{\numvec}[1]{{\vec{\num{#1}}}}
\newcommand{\nummat}[1]{\numvec{#1}}

\newcommand{\boundary}{\partial}

 \newcommand{\CC}{\ensuremath{\operatorname C}\xspace}
 \newcommand{\HH}{\ensuremath{\operatorname H}\xspace}
 \newcommand{\LL}{\ensuremath{\operatorname L}\xspace}
 
 \newcommand{\WW}{\ensuremath{\operatorname W}\xspace}

 \newcommand{\cont}[1]{\ensuremath{\CC^{#1}}}

 \newcommand{\leb}[1]{\ensuremath{\LL_{#1}}}

 \newcommand{\sob}[2]{\ensuremath{\WW^{#1}_{#2}}}
 
 \newcommand{\sobh}[1]{\ensuremath{\HH^{#1}}}
 \newcommand{\sobhz}[1]{\sobh{#1}_0}

 \newcommand{\poly}[1]{\ensuremath{\rP}^{#1}}

 \newcommand{\fes}[1]{\ensuremath{\fespace^{#1}}}

\newcommand{\Forall}{\:\forall\:}

\newcommand{\Foreach}{\quad\Forall}

\usepackage{amscd}

\newcommand{\funk}[3]{\ensuremath{#1:#2\to#3}}

\newcommand{\dfunkmapsto}[6][]{\ensuremath{
    \begin{array}{rccl}
      {#2}: & {#4} &  \to   & {#6}
      \\
            & {#3} &\mapsto & {#5#1}
    \end{array}\quad}}

\renewcommand{\restriction}[2]{\left.#1\right|_{#2}}

\newcommand{\aka}[1]{(also known as {#1})\xspace}

\newcommand{\Program}[1]{\textsf{#1}\xspace}

\newcommand{\matlab}{\Program{Matlab\registered}}

\providecommand{\ListParameters}{}
\renewcommand{\ListParameters}
{
	 \setlength{\topsep}{0em}
	 \setlength{\leftmargin}{0em}
         \setlength{\itemsep}{0ex}
	 \setlength{\parsep}{.5ex}
	 \setlength{\itemindent}{\labelsep}
	 \addtolength{\itemindent}{\labelwidth}
}

\newcounter{LetterListItem}
\renewcommand{\theLetterListItem}{(\alph{LetterListItem})}
{
	\begin{list}%
	{\theLetterListItem\ }%
	{\usecounter{LetterListItem}
	 \ListParameters
	}
}%
{\end{list}}
\newcounter{NumberListItem}
\renewcommand{\theNumberListItem}{\arabic{NumberListItem}}
{
	\begin{list}%
	{\theNumberListItem.\ }%
	{\usecounter{NumberListItem}%
	 \ListParameters
	}
}%
{\end{list}}
\newcounter{QuestionListItem}
\renewcommand{\theQuestionListItem}{\textbf{Question \arabic{QuestionListItem}}}
{
	\begin{list}%
	{\theQuestionListItem.\ }%
	{\usecounter{QuestionListItem}%
	 \ListParameters
	}
}%
{\end{list}}
\newcounter{RomanListItem}
\renewcommand{\theRomanListItem}{(\roman{RomanListItem})}
{
	\begin{list}%
	{\theRomanListItem\ }%
	{\usecounter{RomanListItem}
	 \ListParameters
	}
}%
{\end{list}}
\newcounter{StepsItem}
{
	\begin{list}%
	{Step \theStepsItem.\ }%
	{\usecounter{StepsItem}%
	 \ListParameters
	}
}%
{\end{list}}

\providecommand{\grad}{\nabla}
\renewcommand{\grad}{\nabla}

\usepackage{amsthm} 
\ifthenelse{\boolean{isthesis}}{%
  \setcounter{secnumdepth}{1}%
}{
  \setcounter{secnumdepth}{2} 
}
\providecommand{\ListParameters}{}
\renewcommand{\ListParameters}
{
	 \setlength{\topsep}{0em}
	 \setlength{\leftmargin}{0em}
         \setlength{\itemsep}{0ex}
	 \setlength{\parsep}{.5ex}
	 \setlength{\itemindent}{\labelsep}
	 \addtolength{\itemindent}{\labelwidth}
}
\newtheoremstyle{plain}
  {}
  {}
  {\mdseries\slshape}
  {\parindent}
  {\bfseries}
  {.}
  {.5em}
  {}

\newtheoremstyle{note}
  {}
  {}
  {}
  {\parindent}
  {\bfseries}
  {.}
  {.5em}
  {}

\newtheoremstyle{claim}
  {}
  {}
  {\mdseries\slshape}
  {}
  {\bfseries}
  {:}
  {.5em}
  {}

\newtheoremstyle{exercise}
  {}
  {}
  {}
  {}
  {\bfseries}
  {.}
  {1em}
  {}

\newtheoremstyle{break}
  {}
  {}
  {}
  {}
  {\bfseries}
  {.}
  {\newline}
  {}

  \newcommand{\Proofname}{Proof}

\ifthenelse{\boolean{isthesis}}{%
  
}{%
  
}

\newcommand{\pdfformat}[1]{
  \provideboolean{pdfoutput}
  \setboolean{pdfoutput}{#1}
  \ifthenelse{\boolean{pdfoutput}}%
	     {\typeout{using pdf}
	       \usepackage[pdftex]{graphicx,xcolor}
	       \newcommand{\graphext}{pdf}
	       \newcommand{\graphextex}{pdf_t}
	       \usepackage{epsfig}
	       \usepackage{tikz}
	     }
	     {
	       \typeout{using eps}
	       \usepackage[dvips]{graphicx,xcolor}
	       \newcommand{\graphext}{eps}
	       \newcommand{\graphextex}{eps_t}
	       \usepackage{epsfig}
	       \usepackage{tikz}
	     }
}

\newcommand{\Le}[1]{\ensuremath{\leb{#1}}} 
\newcommand{\hoz}{\sobhz1(\W)}

\newcommand{\T}[1]{\cT^{#1}}

\newcommand{\linop}{\ensuremath{\cL}}

\newcommand{\tr}[1]{\ensuremath{\operatorname{trace}{(#1)}}}

\renewcommand{\rangefromto}[3]{\ensuremath{#1=#2,\ldots,#3}}
 \newcommand{\figscale}{1}

\renewcommand{\fes}{\ensuremath{\mathbb V}}
\renewcommand{\sin}[1]{\ensuremath{\operatorname{sin}\left(#1\right)}}

\renewcommand{\arctan}[1]{\ensuremath{\operatorname{arctan}\left(#1\right)}}

\newcommand{\con}[1]{\ensuremath{\kappa(#1)}}

\newcommand{\tensorp}[2]{\ensuremath{\left\langle#1\otimes#2\right\rangle}}
\newcommand{\frob}[2]{\ensuremath{{#1}{:}{#2}}}
\renewcommand{\div}[1]{\operatorname{div}\left(#1\right)}
\renewcommand{\H}{\ensuremath{\vec{H}}}

\newcommand{\symm}{\ensuremath{\operatorname{Sym}(\reals^{d\times d})}}

\newcommand{\feszero}{\mathring{\fes}}
\newcommand{\Nzero}{\mathring{N}}
\newcommand{\Nd}{{N}_{\pd{}{}}}

\newcommand{\Phizero}{\mathring{\Phi}}

\newcommand{\A}{\ensuremath{\vec{A}}}

\numberwithin{equation}{section}
\setboolean{showtodo}{true}
\setboolean{shownotes}{true}
\setboolean{showchanges}{true}
\setboolean{usemathrsfs}{false}
\author{Omar Lakkis}
\address{ Omar Lakkis\newline
Department of Mathematics\newline
University of Sussex\newline
Brighton\newline
UK-BN1 9RF, United Kingdom} \curraddr{}
\email{\linkedemail{o.lakkis@sussex.ac.uk}}
\urladdr{\linkedurl{http://www.maths.sussex.ac.uk/Staff/OL}}
\author{Tristan Pryer} 
\address{ Tristan Pryer\newline
  Department of Mathematics\newline 
  University of Sussex\newline
  Brighton\newline 
  UK-BN1 9RF, United Kingdom} 
\curraddr{}
\email{\linkedemail{tmp20@sussex.ac.uk}}
\urladdr{\linkedurl{http://www.maths.sussex.ac.uk/~tristan}}
\title[A FEM for nonvariational elliptic problems]
      {A finite element method for second order nonvariational elliptic problems} 
\date{\today}
\pdfformat{true}
\begin{document}
\maketitle
\begin{abstract}
  We propose a numerical method to approximate the solution of second
  order elliptic problems in nonvariational form.  The method is of
  Galerkin type using conforming finite elements and applied directly
  to the nonvariational (nondivergence) form of a second order linear
  elliptic problem.  The key tools are an appropriate concept of
  ``finite element Hessian'' and a Schur complement approach to
  solving the resulting linear algebra problem. The method is
  illustrated with computational experiments on three linear and one
  quasilinear PDE, all in nonvariational form.
\end{abstract}
\section{Introduction}
\label{sec:intro}
  Finite element methods (FEM) arguably constitute one of the most
  successful method families in numerically approximating elliptic
  partial differential equations (PDE's) that are given in variational
  \aka{divergence} form.

  For the reader's appreciation of this statement we briefly introduce
  standard FEM concepts. Let $\W$ be a given domain (open and bounded
  set) in $\R d$, $d\in\naturals$,
  $\funk{f,a_{\a,\b}=a_{\b,\a}}\W\reals$, be given functions with the
  appropriate regularity such that the operator $\div{\A\grad u}$, for
  $\A:=[a_{\a,\b}]_{\rangefromto{\a,\b}1d }$, makes sense, is elliptic
  and there is a unique function $\funk u\W\reals$ satisfying
  $\div{\A\grad u}=f$ with $u=0$ on $\partial\W$~\cite[for
    details]{Gilbarg:1983}.  The \emph{classical solution}, $u$, of
  this problem can be characterized by first writing the PDE in
  \emph{weak} \aka{\emph{variational}} \emph{form} using Green's
  formula:
  \begin{equation}
    u\in\cY\text{ and satisfies } 
    a(u,v)
    :=\int_\W{\Transpose{\grad u}\A\grad v}
    =\int_\W fv\Foreach v\in\cX,
  \end{equation}
  where $\cX \AND \cY$ are appropriate (infinite dimensional) function
  spaces. A (finite) Galerkin procedure consists in finding an
  \emph{approximation} of $u$, $U\in\rY$
  \begin{equation}
    A(U,V)=\ltwop fV\Foreach V\in\rX,
  \end{equation}
  where $\rY$ and $\rX$ are finite dimensional ``counterparts''
  (usually subspaces, but may be not) of $\cY$ and $\cX$ and the
  bilinear form $A$ an approximation of $a$.  For example, when $a=A$
  (modulo quadrature) $\cX=\cY=\sobhz1(\W)$ and $\rX=\rY$ are a space
  of continuous piecewise $p$-degree polynomial functions on a
  partition of $\W$, we obtain the standard \emph{conforming
   mesh-refinement ($h$-version) finite element method of degree $p$}.

  The reason behind the FEM's success in such a framework is twofold:
  (1) the weak form is suitable to apply functional analytic
  frameworks (Lax--Milgram Theorem or Babu{\v{s}}ka--Brezzi--Lady{\v z}enskaya
  condition, e.g.), and (2) the discrete functions need to be
  differentiated at most once, whence weak smoothness requirements on
  the ``elements''.

  In this article, we depart from this basis by considering second
  order elliptic boundary value problems (BVP's) in nonvariational
  form
  \begin{equation}
    \label{eqn:BVP:strong-form}
    \text{find }u\text{ such that }
    \frob{\A}{\Hess u} = f\text{ in }\W
    \AND
    \restriction u{\boundary\W}=g,
  \end{equation}
  for which one may not always be successful in applying the standard
  FEM (with reference to \S\ref{sec:notation} for the notation).
  Indeed, the use of the standard FEM requires (1) the coefficient
  matrix $\funk \A\W{\R{d\times d}}$ to be (weakly) differentiable and
  (2) the rewriting of the second order term in divergence form, an
  operation which introduces an advection (first order) term:
  \begin{equation}
    \label{eq:divform}
    \frob{\A}{\Hess u} = \div{\A \nabla u} - \qp{\div \A} \nabla u.
  \end{equation}
  Even when coefficient matrix $\mat A$ is differentiable on $\W$,
  this procedure could result in the problem becoming
  advection--dominated and unstable for conforming FEM, as we
  demonstrate numerically using Problem (\ref{eq:Problem2}).
  
  Our main motivation for studying linear elliptic BVP's in
  nonvariational form is their important role in pure and applied
  mathematics. An important example of nonvariational problems is the
  fully nonlinear BVP that is approximated via a Newton method which
  becomes an infinite sequence of linear nonvariational elliptic
  problems~\cite{Bohmer:2008}.
  
  In this article, we propose and test a direct discretization of the
  \emph{strong form} (\ref{eqn:BVP:strong-form}) that makes no special
  assumption on the derivative of $\A$. The main idea, is an
  appropriate definition of a \emph{finite element Hessian} given in
  \S\ref{defn:fehessian}. The finite element Hessian has been used
  earlier in different contexts, such as anisotropic mesh generation
  \cite{AgouzalVassilevski:2002,ChenSunXu:2007,ValletManoleDompierreDufourGuibault:2007}
  and \emph{finite element convexity}~\cite{Aguilera:2008}.  The
  finite element Hessian is related also to the finite element
  (discrete) elliptic operator appearing in the analysis of evolution
  problems~\cite{Thomee:2006}.
  
  The method we propose is quite straightforward, and we are surprised
  that it is not easily available in the literature.  It consists in
  discretizing, via a Galerkin procedure, the
  BVP~(\ref{eqn:BVP:strong-form}) \emph{directly without writing it
  in divergence form}.

  The main difficulty of our approach is having to deal with a
  somewhat involved linear algebra problem that needs to be solved as
  efficiently as possible (this is especially important when we apply
  this method in the linearization of nonlinear elliptic BVP's).  We
  overcame this difficulty in \S\ref{sec:characteriseandsolve}, by combining
  the definition of $u$'s distributional Hessian,
  \begin{equation}
    \duality{\Hess u}{\phi} 
    =
    - \tensorp{\nabla u}{\nabla\phi}
    + \tensorp{\nabla u}{\vec n \ \phi}_{\pd{}\W}    
    \Foreach \phi\in\cont{\infty}(\W),
  \end{equation}
  with equation (\ref{eqn:BVP:strong-form}) into a system of equations
  that are larger, but easier to handle numerically, once discretized.

  It is worth noting that there are alternatives to our approach, most
  notably the standard finite difference method and its variants.  The
  reason we are interested in a Galerkin procedure is the ability to
  use an unstructured mesh, essential for complicated geometries where
  the finite difference method leads to complicated, and sometimes
  prohibitive, modifications (especially in dimension $3$ and higher),
  and the potential of dealing with adaptive methods, using available
  finite element code.  Furthermore, our method has the potential to
  approach the iterative solution fully nonlinear problems where
  finite difference methods can become clumsy and
  demanding~\cite{KuoTrudinger:1992, LoeperRapetti:2005, Oberman:2008,
    CaffarelliSouganidis:2008}.

  This paper focuses mainly on the algorithmic and linear algebraic
  aspects of the method and is set out as follows. In
  \S\ref{sec:notation} we introduce some notation and set out the
  model problem. We then present a discretization scheme for the model
  problem using standard conforming finite elements in
  $\cont{0}(\W)$. In \S\ref{sec:characteriseandsolve} we present a
  linear algebra technique, inspired by the standard \emph{Schur
    complement} idea, for solving the linear system arising from the
  discretization. Finally, in \S\ref{sec:numerics} we summarize
  extensive numerical experiments on model linear boundary value
  problems (BVPs) in nonvariational form and an application to
  quasilinear BVP in nonvariational form.
  \clearpage
\section{Set up}
\label{sec:notation}
\subsection{Notation}
\label{sse:notation}
Let $\W\subset \R{d}$ be an open and bounded Lipschitz domain. We
denote $\leb{2}(\W)$ to be the space of square (Lebesgue) integrable
functions on $\W$ together with it's inner product $\ltwop{v}{w} :=
\int_\W v w$ and norm $\Norm{v} := \Norm{v}_{\leb2(\W)} =
\ltwop{v}{v}^{1/2}$. We also denote by $\langle f \rangle_\w$ the
integral of a function $f$ over the domain $\w$ and drop the subscript
for $\w = \W$.

We use the convention that the derivative $\D u$ of a function
$u:\W\to\reals$ is a row vector, while the gradient of $u$, $\nabla u$
is the derivative's transpose, i.e., $\nabla u = \Transpose{\left(\D
  u\right)}$. We will make use of the slight abuse of notation,
following a common practice, whereby the Hessian of $u$ is denoted as
$\Hess u$ (instead of the correct $\nabla \D u$) and is represented by
a $d\times d$ matrix.

The Sobolev spaces
~\cite{Ciarlet:1978,Evans:1998}
\begin{equation}
  \sobh{k}(\W)
  := 
  \sob{k}{2}(\W) 
  = 
  \ensemble{\phi\in\leb2(\W)}
           {\sum_{\norm{\vec\alpha} \leq k} \D^{\vec\alpha}\phi\in\leb2(\W)},
\end{equation}
are equipped with norms and semi-norms
\begin{gather}
    \Norm{v}_{k}^2
    := 
    \Norm{v}_{\sobh{k}(\W)}^2 
    = 
    \sum_{\norm{\vec \alpha}\leq k}\Norm{\D^{\vec \alpha} v}^2 
    \\
    \AND \norm{v}_k^2 
    :=
    \norm{v}_{\sobh{k}(\W)}^2 
    =
    \sum_{\norm{\vec \alpha} = k}\Norm{\D^{\vec \alpha} v}^2
\end{gather}
respectively, where $\vec\alpha = \{ \alpha_1,...,\alpha_d\}$ is a
multi-index, $\norm{\vec\alpha} = \sum_{i=1}^d\alpha_i$ and
derivatives $\D^{\vec\alpha}$ are understood in a weak sense. We pay
particular attention to the cases $k = 1,2$,
\begin{gather}
  \hoz := \text{closure of }\cont{\infty}_0(\W) \text{ in } \sobh{1}(\W)
  \\
  \AND \sobh{-1}(\W):=\operatorname{dual}\qp{\hoz}.
\end{gather}
We denote by $\duality{v}{w}$ the action of a distribution $v$ on the
function $w$. If both $v,w\in\Le{2}(\W)$ then $\duality{v}{w} =
\ltwop{v}{w}$.

We consider the following problem: Find $u\in\hoz$ such that
\begin{equation}
  \begin{split}
    \label{Problem}
    \linop u &= f \text{ in }\W,\\
    u &= 0 \text { on }\partial \W,
  \end{split}
\end{equation}
where the data $f : \W \to \R{}$ is prescribed and $\linop$ is a
general linear, second order, uniformly elliptic partial differential
operator. Let $\A \in \Le{\infty}(\W)^{d\times d} \cap \symm$, the
space of bounded, symmetric, positive definite, $d \times d$ matrixes.
\begin{equation}
  \dfunkmapsto[,]
	      {\linop}
	      u
	      {\hoz}
	      {\linop u:= \frob{\A}{\Hess u}}
	      {\sobh{-1}(\W)}
\end{equation}
we use $\frob{\vec X}{\vec Y}:=\tr{\Transpose{\vec X} \vec Y}$ to
denote the Frobenius inner product between two matrixes.
\subsection{Discretization}
\label{sec:discretisation}
Let $\T{}$ be a conforming triangulation of $\W$, namely, $\T{}$ is a
finite family of sets such that
\begin{enumerate}
\item $K\in\T{}$ implies $K$ is an open simplex (segment for $d=1$,
  triangle for $d=2$, tetrahedron for $d=3$),
\item for any $K,J\in\T{}$ we have that $\closure K\meet\closure J$ is
  a full subsimplex (i.e., it is either $\emptyset$, a vertex, an
  edge, a face, or the whole of $\closure K$ and $\closure J$) of both
  $\closure K$ and $\closure J$ and
\item $\union{K\in\T{}}\closure K=\closure\W$.
\end{enumerate}
The \emph{shape regularity} of $\T{}$ is defined as
\begin{equation}
  \label{eqn:def:shape-regularity}
  \mu(\T{}) := \inf_{K\in\T{}} \frac{\rho_K}{h_K},
\end{equation}
where $\rho_K$ is the radius of the largest ball contained inside
$K$ and $h_K$ is the diameter of $K$.
 We use the convention where $\funk h\W\reals$ denotes the
 \emph{meshsize function} of $\T{}$, i.e.,
 \begin{equation}
   h(\vec{x}):=\max_{\closure K\ni \vec x}h_K.
 \end{equation}
 We introduce the \emph{finite element spaces}
 \begin{gather}
   \label{eqn:def:finite-element-space}
   \fes
   :=\ensemble{\Phi \in \sobh1(\W)}{\Phi\vert_{K} \in \poly p\Forall K\in\T{}},
   \\
   \feszero
   :=\fes \cap \hoz,
 \end{gather}
 where $\poly k$ denotes the linear space of polynomials in $d$
 variables of degree no higher than a positive integer $k$. We
 consider $p\geq 1$ to be fixed and denote by $\Nzero :=
 \dim{\feszero}$ and $N = \Nzero + \Nd := \dim{\fes}$. Let $\numvec{\Phizero} =
 \Transpose{( \Phizero_1, ..., \Phizero_{\Nzero} )}$ and $\numvec \Phi
 = \Transpose{( \Phizero_1, ..., \Phizero_{\Nzero}, \Phi_1, ...,
   \Phi_{\Nd} )}$ where $ \{ \Phizero_1, ..., \Phizero_{\Nzero}\}$ and $
 \{\Phizero_1, ..., \Phizero_{\Nzero}, \Phi_1, ...,
   \Phi_{\Nd} \}$  form a basis of $\feszero$, $\fes$ 
 respectively.

Testing the model problem (\ref{Problem}) with
$\phi\in\hoz$ gives
\begin{equation}
  \label{FullProblem}
  \ltwop{\linop u}{\phi} = \ltwop{\frob{\A}{\Hess u}}{\phi} = \ltwop{f}{\phi}.
\end{equation}
In order to discretize \eqref{FullProblem} with $\fes$ we use
an appropriate definition of a Hessian of a finite element
function. Such a function may not admit a Hessian in the classical
sense, so we consider it as a distribution (or generalized function)
which we recall the definition.
\begin{Defn}[generalized Hessian]
  \label{defn:distributionalhessian}
  Let $\funk{\vec n}{\boundary\W}{\R d}$ be the outward pointing normal of $\W$. Given
  $v\in\hoz$ its \emph{generalized Hessian} defined in the standard
  distributional sense is given by
  \begin{equation}
      \duality{\Hess v}{\phi} 
      =
      - \tensorp{\nabla v}{\nabla\phi}
      + \tensorp{\nabla v}{\vec n \ \phi}_{\boundary\W}    
      \Foreach \phi\in\cont{\infty}(\W),
  \end{equation}
  where  we
  are using $\vec x \otimes \vec y:={\vec x}{\Transpose{\vec y}}$ to
  denote the tensor product between two geometric vectors $\vec x \AND
  \vec y$.
\end{Defn}

\begin{The}[finite element Hessian]
  \label{the:fehessian}
  For each $V\in\feszero$ there exists a unique
  $\H[V]\in\fes^{d\times d}$ such that
  \begin{equation}
    \ltwop{\H[V]}{\Phi} 
    = 
    \duality{\Hess V}{\Phi} \Foreach \Phi\in\fes.
  \end{equation}
\end{The}

\begin{Proof}
  Given a finite element function $V\in\feszero$, Definition
  \ref{defn:distributionalhessian} implies
  \begin{equation}
    \duality{\Hess V}{\phi} 
    = 
    - \tensorp{\nabla V}{\nabla\phi}
    + \tensorp{\nabla V}{\vec n \ \phi}_{\boundary\W} 
    \Foreach \phi\in\cont{\infty}(\W).
  \end{equation}
  We fix $V$ and let
  \begin{equation}
    \dfunkmapsto[.]
                {G}
                {\phi}
                {\cont{\infty}(\W)}
                {-\tensorp{\nabla V}{\nabla \phi}    
                  + \tensorp{\nabla V}{\vec n \ \phi}_{\boundary\W} 
                }
                {\reals^{d\times d}}
  \end{equation}
  Notice that $G$ is a bounded linear functional on
  $\cont{\infty}(\W)$ in the $\sobh1(\W)$-norm as,
  \begin{equation}
      \norm{G(\phi)} 
      =
      \norm{\tensorp{\nabla V}{\nabla \phi}}
      + 
      \norm{\tensorp{\nabla V}{\vec n \ \phi}_{\boundary\W} }
      \leq
      C(d,\W) \Norm{V}_1\Norm{\phi}_1.
  \end{equation}
  Thus, due to the density of $\cont{\infty}(\W)$ in $\sobh1(\W)$, $G$
  admits a unique extension, $\tilde G$.

  Let $R = \restriction{\tilde G}{\fes}$ be the restriction of $\tilde
  G$ to $\fes$. Since $\tilde G$ is linear and bounded on $\sobh1(\W)$
  it follows that $R$ is linear and bounded on $\fes$ in the
  $\sobh1(\W)$-norm. Hence by Riesz's Representation Theorem there
  exists an $\H[V]\in\fes^{d\times d}$ such that for each
  $\Phi\in\fes$
  \begin{equation}
      \ltwop{\H[V]}{\Phi} 
      :=
      R(\Phi)
      = 
      -\tensorp{\nabla V}{\nabla \Phi} 
      +
      \tensorp{\nabla V}{\vec n \ \Phi}_{\boundary\W} 
      ,
  \end{equation}
  which coincides with the generalized Hessian (cf.
  Definition \ref{defn:distributionalhessian}) on $\fes$.
\end{Proof}

\begin{Defn}[finite element Hessian]
  \label{defn:fehessian}
  From Theorem \ref{the:fehessian} we define the
  \emph{finite element Hessian} as follows. Let $V\in\feszero$ then
  \begin{equation}
    \ltwop{\H[V]}{\Phi} 
    :=
    - \tensorp{\nabla V}{\nabla\Phi}
    +
    \tensorp{\nabla V}{\vec n \ \Phi}_{\boundary\W} 
    \Foreach \Phi\in\fes.
  \end{equation}
  It follows that $\H$ is a linear
  operator on $\feszero$.
\end{Defn}
Taking the model problem (\ref{FullProblem}) we substitute the
finite element Hessian directly, reducing the space of test
functions to $\feszero$, we wish to find $U\in\feszero$ such that
\begin{equation}
  \label{eqn:disc}
  \ltwop{\frob{\A}{\H[U]}}{\Phizero} 
  =
  \ltwop{f}{\Phizero} \Foreach \Phizero\in\feszero.
\end{equation}
\begin{The}[nonvariational finite element method (NVFEM)]
  \label{the:ndfem}
    The \emph{nonvariational finite element solution} for the model
    problem's discretization (\ref{eqn:disc}) is given as $U =
    \Transpose{\numvec \Phizero} \numvec u$, where $\numvec u \in
    \reals^{\Nzero}$ is the solution to the following linear system
  \begin{equation}
    \label{eqn:ndfemformulation}
    \nummat D \numvec u 
    := 
    \sum_{\alpha=1}^d\sum_{\beta=1}^d \nummat B^{\alpha,\beta} \nummat
    M^{-1} \nummat C_{\alpha,\beta} \numvec{u} = \numvec{f}.
  \end{equation}
  The components of (\ref{eqn:ndfemformulation}) are given by 
  \begin{align}
    \label{def:probmat}
    \nummat B^{\alpha,\beta} 
    &:=
    \ltwop{\numvec \Phizero}{\A^{\alpha,\beta}\Transpose{\numvec \Phi}}
    \in \reals^{\Nzero \times N},
    \\
    \label{def:mass}
    \nummat{M}
    &:=
    \ltwop{\numvec \Phi}{\Transpose{\numvec\Phi}}
    \in \reals^{N\times N},
    \\
    \label{def:stiff}
    \nummat C_{\alpha,\beta}
    &:=
    -
    \ltwop{\pd{\beta}{\numvec\Phi} }
         {\pd\alpha{}\Transpose{\numvec\Phizero}}
    + 
    \ltwop{\numvec\Phi n_\beta}
          {\pd\alpha{\Transpose{\numvec\Phizero}}}_{\boundary\W}
          \in \R{N\times \Nzero},
    \\
    \numvec f
    &:= 
    \ltwop{f}{{\numvec\Phizero}}
    \in \reals^{\Nzero}.
  \end{align}
\end{The}

\begin{Proof}
  Since $\H[U] \in \fes^{d\times d}$ for each $\alpha,\beta =
  1,\dots,d$ , $\H_{\alpha,\beta}[U] = \Transpose{\numvec \Phi }
  \numvec h_{\alpha,\beta} $. Then, testing (\ref{eqn:disc}) with
  $\numvec \Phizero$,
\begin{equation}
  \label{der:hessian}
  \begin{split}
    \ltwop{f}{\numvec \Phizero} 
    &= 
    \sum_{\alpha=1}^d\sum_{\beta=1}^d
    \ltwop{\A^{\alpha,\beta}\H_{\alpha,\beta}[U]} {\numvec \Phizero}
    \\
    &=
    \sum_{\alpha=1}^d\sum_{\beta=1}^d 
    \ltwop{\numvec \Phizero } {\A^{\alpha,\beta} \Transpose{\numvec \Phi}
      \numvec
      h_{\alpha,\beta} }
    \\ 
    &=
    \sum_{\alpha=1}^d\sum_{\beta=1}^d 
    \ltwop{ \numvec \Phizero}{\A^{\alpha,\beta} \Transpose{\numvec \Phi}}
    {\numvec h_{\alpha,\beta}}.
    \\
    &= 
    \sum_{\alpha=1}^d\sum_{\beta=1}^d 
    \nummat B^{\alpha,\beta} \numvec h_{\alpha,\beta}
  \end{split}
\end{equation}

Utilizing Definition \ref{defn:fehessian} for each $\alpha, \beta =
1\dots d$ we can compute $\numvec h_{\alpha,\beta} \in\reals^N$,
noting $U = \Transpose{\numvec{\Phizero}}\numvec{u}$,
\begin{equation}
  \begin{split}
    \ltwop{\numvec\Phi}{\Transpose{\numvec\Phi}}
    \numvec h_{\alpha,\beta}
    &=
    \ltwop{\numvec\Phi}{\H_{\alpha,\beta}[U]}
    \\
    &=
    -
    \ltwop{\pd{\beta}{}\numvec\Phi }{\pd\alpha U} 
    + 
    \ltwop{\numvec\Phi \vec n_\beta}{\pd\alpha U}_{\boundary\W}
    \\
    &=
    \left(-\ltwop{\pd{\beta}{}\numvec\Phi }
         {\pd\alpha{\Transpose{\numvec\Phizero}}}
    + 
    \ltwop{\numvec\Phi \vec n_\beta}
          {\pd\alpha{\Transpose{\numvec\Phizero}}}_{\boundary\W}\right)
          \numvec u.
  \end{split}
\end{equation}
Using the definition of $\nummat C_{\alpha,\beta}$ (\ref{def:stiff})
and $\nummat M$ (\ref{def:mass}) we see for each $\alpha, \beta =
1\dots d$
\begin{equation}
  \label{der:hessiancoeff}
  \begin{split}
    \nummat M  {\numvec h_{\alpha,\beta}}
    &= 
    \nummat C_{\alpha,\beta} {\numvec{u}}
    \\
    {\numvec h_{\alpha,\beta}} 
    &=
    \nummat M^{-1}\nummat C_{\alpha,\beta}{\numvec{u}}  .
  \end{split}
\end{equation}
Substituting $\numvec h_{\alpha,\beta}$ from \eqref{der:hessiancoeff}
into \eqref{der:hessian} we obtain the desired result.
\end{Proof}

\begin{Example}[for $d=2$]
  For a general elliptic operator in 2-D, the formulation
  \eqref{eqn:ndfemformulation} takes the form
  \begin{equation}
    \left( \nummat B^{1,1} \nummat M^{-1} \nummat C_{1,1} 
    +
    \nummat B^{2,2} \nummat M^{-1} \nummat C_{2,2} 
    + 
    \nummat B^{1,2} \nummat M^{-1} \nummat C_{1,2}
    +
    \nummat B^{2,1} \nummat M^{-1} \nummat C_{2,1}
    \right) \numvec{u}
    =  \numvec{f}
  \end{equation}
\end{Example}

\section{Solving the linear system}
\label{sec:characteriseandsolve}

\begin{Obs}[(\ref{eqn:ndfemformulation}) is difficult to solve]
  Looking at the full system setting $\nummat D = \sum\sum\nummat
  B^{\alpha,\beta} \nummat M^{-1} \nummat C_{\alpha,\beta}$
  multiplying out each of the matrixes and proceeding to solve
  $\nummat D \numvec u = \numvec f$ the resulting system would not be
  sparse forcing the use of direct solvers.
\end{Obs}

In this section we will present a method to solve formulation
(\ref{eqn:ndfemformulation}) in a general setting. This method makes
use of the sparsity of the component matrixes $\nummat
B^{\alpha,\beta}, \nummat C^{\alpha,\beta}$ and $\nummat M$.

\begin{Rem}
  An interesting point of note is that if the mass matrix $\nummat M$
  were diagonalized, by mass lumping, then for each $\alpha$ and
  $\beta$ the matrix $\nummat B^{\alpha,\beta} \nummat M^{-1} \nummat
  C_{\alpha,\beta}$ would still be sparse (albeit less so than the
  individual matrixes $\nummat B^{\alpha,\beta}$ and $\nummat
  C_{\alpha,\beta}$). Hence the system can be easily solved using
  existing sparse methods. However mass lumping is only applicable to
  $\poly{1}$ finite elements. For higher order finite elements it
  would be desirable to exploit the sparse structure of the component
  matrixes that make up the system.
\end{Rem}

\subsection{A generalized Schur complement}

We observe the matrix $\nummat D$ in the system
(\ref{eqn:ndfemformulation}) is a sum of Schur complements $\nummat
B^{\alpha,\beta}\nummat{M}^{-1}\nummat C_{\alpha,\beta}$. With that in mind we
introduce the $(d^2 + 1)^2$ block matrix
\begin{equation}
  \nummat E =     \left[ 
    \begin{array}{cccccc}
      \nummat M & \nummat 0 & \dotsb & \nummat 0 & \nummat 0 & -\nummat C_{1,1}\\
      \nummat 0 & \nummat M & \dotsb & \nummat 0 & \nummat 0 & -\nummat C_{1,2}\\
      \vdots & \vdots & \ddots & \vdots & \vdots & \vdots\\
      \nummat 0 & \nummat 0 & \dotsb & \nummat M & \nummat 0 & -\nummat C_{d,d-1}\\
      \nummat 0 & \nummat 0 & \dots & \nummat 0 & \nummat M & -\nummat C_{d,d}\\
      \nummat B^{1,1} & \nummat B^{1,2}  & \dots & \nummat B^{d,d-1} & \nummat B^{d,d} & \nummat 0
    \end{array} 
    \right].
\end{equation}
\begin{Lem}[generalized Schur complement]
  Given 
  \begin{gather}
    \numvec v = \Transpose{
      \left(
      \numvec h_{1,1}, \numvec h_{1,2},
      \dots, \numvec h_{d,d-1}, \numvec h_{d,d}, \numvec u
      \right)},
    \\
    \numvec b = \Transpose{
      \left(
      \numvec 0, \numvec 0
      \dots,
      \numvec 0, \numvec 0, \numvec f
      \right)},
  \end{gather}
  solving the system 
  \begin{equation}
    \nummat D \numvec u = 
    \sum_{\alpha=1}^d\sum_{\beta=1}^d \nummat B^{\alpha,\beta} \nummat
    M^{-1} \nummat C_{\alpha,\beta} \numvec{u} 
    =
    \numvec{f},
  \end{equation}
  is equivalent to solving
  \begin{equation}
    \label{eq:blocksystem}
    \nummat E \numvec v = \numvec b.
  \end{equation}
  for $\numvec u$.
\end{Lem}
\begin{Proof}
  The proof is just block Gaussian elimination on $\nummat
  E$. Left-multiplying the first $d^2$ rows by $\nummat M^{-1}$ yields
  \begin{equation}
    \left[ 
      \begin{array}{cccccc}
        \nummat I & \nummat 0 & \dotsb & \nummat 0 & \nummat 0 & -\nummat M^{-1}\nummat C_{1,1}\\
        \nummat 0 & \nummat I & \dotsb & \nummat 0 & \nummat 0 & -\nummat M^{-1}\nummat C_{1,2}\\
        \vdots & \vdots & \ddots & \vdots & \vdots & \vdots\\
      \nummat 0 & \nummat 0 & \dotsb & \nummat I & \nummat 0 & -\nummat M^{-1}\nummat C_{d,d-1}\\
      \nummat 0 & \nummat 0 & \dots & \nummat 0 & \nummat I & -\nummat M^{-1}\nummat C_{d,d}\\
      \nummat B^{1,1} & \nummat B^{1,2}  & \dots & \nummat B^{d,d-1} & \nummat B^{d,d} & \nummat 0
      \end{array} 
    \right]
    \left[
      \begin{array}{c}
        \numvec h_{1,1}\\
        \numvec h_{1,2}\\
        \vdots\\
        \numvec h_{d,d-1}\\
        \numvec h_{d,d}\\
        \numvec u
      \end{array}
      \right]
    =
    \left[
      \begin{array}{c}
        \numvec 0\\
        \numvec 0\\
        \vdots\\
        \numvec 0\\
        \numvec 0\\
        \numvec f
      \end{array}
      \right].
    \end{equation}
    Multiplying the $i$-th row by the $i$-th entry of the $(d^2+1)$-th row
    for $i = 1,\dots, d^2$
    \begin{equation}
      \left[ 
        \begin{array}{cccccc}
        \nummat  B^{1,1} & \nummat 0 & \dotsb & \nummat 0 & \nummat 0 & -\nummat  B^{1,1}\nummat M^{-1}\nummat C_{1,1}\\
        \nummat 0 & \nummat B^{1,2} & \dotsb & \nummat 0 & \nummat 0 & -\nummat B^{1,2}\nummat M^{-1}\nummat C_{1,2}\\
        \vdots & \vdots & \ddots & \vdots & \vdots & \vdots\\
      \nummat 0 & \nummat 0 & \dotsb & \nummat B^{d,d-1} & \nummat 0 & -\nummat B^{d,d-1} \nummat M^{-1}\nummat C_{d,d-1}\\
      \nummat 0 & \nummat 0 & \dots & \nummat 0 & \nummat B^{d,d} & -\nummat B^{d,d}\nummat M^{-1}\nummat C_{d,d}\\
      \nummat B^{1,1} & \nummat B^{1,2}  & \dots & \nummat B^{d,d-1} & \nummat B^{d,d} & \nummat 0
      \end{array} 
    \right]
    \left[
      \begin{array}{c}
        \numvec h_{1,1}\\
        \numvec h_{1,2}\\
        \vdots\\
        \numvec h_{d,d-1}\\
        \numvec h_{d,d}\\
        \numvec u
      \end{array}
      \right]
    =
    \left[
      \begin{array}{c}
        \numvec 0\\
        \numvec 0\\
        \vdots\\
        \numvec 0\\
        \numvec 0\\
        \numvec f
      \end{array}
      \right].
    \end{equation}
    Subtracting each of the first $d^2$ rows from the $(d^2+1)$-th row
    reduces the system into row echelon form.
    \begin{equation}
      \left[ 
        \begin{array}{cccccc}
        \nummat  B^{1,1} & \nummat 0 & \dotsb & \nummat 0 & \nummat 0 & -\nummat  B^{1,1}\nummat M^{-1}\nummat C_{1,1}\\
        \nummat 0 & \nummat B^{1,2} & \dotsb & \nummat 0 & \nummat 0 & -\nummat B^{1,2}\nummat M^{-1}\nummat C_{1,2}\\
        \vdots & \vdots & \ddots & \vdots & \vdots & \vdots\\
      \nummat 0 & \nummat 0 & \dotsb & \nummat B^{d,d-1} & \nummat 0 & -\nummat B^{d,d-1} \nummat M^{-1}\nummat C_{d,d-1}\\
      \nummat 0 & \nummat 0 & \dots & \nummat 0 & \nummat B^{d,d} & -\nummat B^{d,d}\nummat M^{-1}\nummat C_{d,d}\\
      \nummat 0 & \nummat 0  & \dots & \nummat 0 & \nummat 0 & \nummat D
      \end{array} 
    \right]
    \left[
      \begin{array}{c}
        \numvec h_{1,1}\\
        \numvec h_{1,2}\\
        \vdots\\
        \numvec h_{d,d-1}\\
        \numvec h_{d,d}\\
        \numvec u
      \end{array}
      \right]
    =
    \left[
      \begin{array}{c}
        \numvec 0\\
        \numvec 0\\
        \vdots\\
        \numvec 0\\
        \numvec 0\\
        \numvec f
      \end{array}
      \right].
    \end{equation}
\end{Proof}

\begin{Rem}[structure of the block matrix]
  In fact this method for the solution of the system $\nummat D
  \numvec u = \numvec f$ is not surprising given the discretization
  presented in the proof of Theorem \ref{the:ndfem} is equivalent to
  the following system:
  \begin{equation}
  \text{Find $U\in\feszero$ such that }
  \begin{cases}
    \ltwop{\H[U]}{\Phi} 
    =
    - \tensorp{\nabla U}{\nabla\Phi}
    +
    \tensorp{\nabla U}{\vec n \ \Phi}_{\pd{}\W} 
    \Foreach \Phi\in\fes
    \\
    \\
    \ltwop{\frob{\A}{\H[U]}}{\Phizero} 
    =
    \ltwop{f}{\Phizero} \Foreach \Phizero\in\feszero.
    \end{cases}
  \end{equation}
\end{Rem}

\begin{Rem}[enforcing non-trivial Dirichlet boundary values]
  \label{rem:non-trivial-dirichlet}
   Given additional problem data $g\in\sobh{1/2}(\W)$, to solve
  \begin{equation}
    \begin{split}
      \label{Problem-nonhomdirichlet}
      \linop u =& f \text{ in }\W,\\
      u =& g \text { on }\partial \W,
    \end{split}
  \end{equation}
  it is not immediate how to enforce the boundary conditions. If we
  were solving the full system $\nummat D \numvec u = \numvec f$, we
  could directly enforce them into the system matrix.

  Since $g\in\sobh{1/2}(\W)$ by an embedding it is
  continuous and can be approximated by the Lagrange interpolant with
  optimal order. To enforce the Dirichlet boundaries we introduce a
  further block representation
  \begin{equation}
    \label{eq:dirichletblock}
    \left[
      \begin{array}{c c}
        \nummat I & \nummat 0
        \\
        \nummat E_{\pd{}{}} & \nummat E
      \end{array}
      \right]
    \left[
      \begin{array}{c}
        \numvec v_{\pd{}{}}
        \\
        \numvec v
      \end{array}
      \right]
    =
    \left[
      \begin{array}{c}
        \numvec b_{\pd{}{}}
        \\
        \numvec b
      \end{array}
      \right],
  \end{equation}
  where $\nummat E, \numvec v \AND \numvec b$ are defined as before
  and $\nummat E_{\pd{}{}}, \numvec v_{\pd{}{}} \AND \numvec
  b_{\pd{}{}}$ are defined as follows
 \begin{gather}
    \nummat E_{\pd{}{}}
    =
    \left[
      \begin{array}{c c c c c c}
        \nummat M & \nummat 0 & \dotsb & \nummat 0 & \nummat 0 & -\nummat C_{1,1}^{\pd{}{}}\\
        \nummat 0 & \nummat M & \dotsb & \nummat 0 & \nummat 0 & -\nummat C_{1,2}^{\pd{}{}}\\
        \vdots & \vdots & \ddots & \vdots & \vdots & \vdots\\
        \nummat 0 & \nummat 0 & \dotsb & \nummat  M & \nummat 0 & -\nummat C_{d,d-1}^{\pd{}{}}\\
        \nummat 0 & \nummat 0 & \dots & \nummat 0 & \nummat M & -\nummat C_{d,d}^{\pd{}{}}\\
        \nummat B^{1,1} & \nummat B^{1,2}  & \dots & \nummat B^{d,d-1} & \nummat B^{d,d} & \nummat 0
      \end{array}
      \right],
    \\
    \numvec v_{\pd{}{}}
    =
    \Transpose{\left[
      \numvec h_{1,1}^{\pd{}{}}, \numvec h_{1,2}^{\pd{}{}}, 
      \dots, \numvec h_{d, d-1}^{\pd{}{}}, \numvec h_{d, d}^{\pd{}{}}, \numvec u^{\pd{}{}}
    \right]},
    \\
    \numvec b_{\pd{}{}}
    =
    \Transpose{\left[
      \nummat 0, \numvec 0, \dots, \numvec 0, \numvec 0, \numvec g
    \right]}.
    \end{gather}
   Let $\nummat \Phi_{\pd{}{}} = \{ \Phi_1,\dots,\Phi_{\Nd} \}$, then
   the components of $\nummat E_{\pd{}{}}$ and $\numvec b_{\pd{}{}}$
   are defined as follows
  \begin{gather}
    \nummat C^{\pd{}{}}_{\alpha,\beta}
    =
    -
    \ltwop{\pd{\beta}{}\numvec\Phi }
         {\pd\alpha{}\Transpose{\numvec \Phi_{\pd{}{}}}}
    + 
    \ltwop{\numvec\Phi \vec n_\beta}
          {\pd\alpha{}\Transpose{\numvec\Phi_{\pd{}{}}}}_{\pd{}\W}
          \in \reals^{N\times \Nd},
          \\
    \numvec g_j 
    =
    g(x_j)\Phi_j
    \in \reals^{\Nd},
  \end{gather}
  where $x_j$ is the Lagrange node associated with $\Phi_j$.

  The block matrix (\ref{eq:dirichletblock}) can then be trivially solved
  \begin{equation}
    \nummat E \numvec v = \numvec b - \nummat E_{\pd{}{}} \numvec b_{\pd{}{}}.
  \end{equation}
\end{Rem}

\begin{Rem}[storage issues]
  We will be using the generalized minimal
  residual method (GMRES) to solve this system. The GMRES, as with any
  iterative solver, only requires an algorithm to compute a
  matrix-vector multiplication. Hence we are only required to store
  the component matrixes $\nummat B^{\alpha,\beta}, \nummat C_{\alpha,
    \beta} \AND \nummat M$.
\end{Rem}

\begin{Rem}[condition number]
  \label{rem:conditionnumber}
   The convergence rate of an iterative solver applied to a linear
   system $\nummat N \numvec v = \numvec g$ will depend on the
   condition number $\con{\nummat N}$, defined as the ratio of the
   maximum and minimum eigenvalues of $\nummat N$:
  \begin{equation}
    \con{\nummat N} = \frac{\lambda_{\text{max}}}{\lambda_{\text{min}}}
  \end{equation}
  Numerically we observe the condition number of the block matrix
  $\con{\nummat E} \leq C h^{-2}$ (see Table \ref{tbl:conditionnumber}).
\end{Rem}

\section{Numerical applications}
\label{sec:numerics}

In this section we study the numerical behavior of the scheme
presented above. All our computations were carried out in \matlab (code
available on request).

We present two linear benchmark problems, for which the solution is known. We
take $\W$ to be the square $S = (-1,1)\times(-1,1) \subset \R 2$ and
in the first two tests consider the operator 
\begin{equation}
  \A(\vec x) =
  \begin{bmatrix}
    1 &  b(\vec x)
    \\
    b(\vec x) & a(\vec x)
  \end{bmatrix}
\end{equation}
varying the coefficients $a(\vec x)$ and $b(\vec x)$. 
\subsection{Test problem with a nondifferentiable operator}
\label{test:nondiff}
For the first test problem we choose the operator in such a way that
\eqref{eq:divform} does not hold, that is the components of $\A$ are
non-differentiable on $\W$, in this case we take
\begin{gather}
  \label{eq:Problem1}
  a(\vec{x})
  = 
  ( x_1^2x_2^2 )^{1/3} + 1
  \\
  b(\vec x) = 0.
\end{gather}
A visualization of the operator \eqref{eq:Problem1} is given in Figure
\ref{Fig:x2y2}. We choose our problem data $f$
such that the exact solution to the problem is given by:
\begin{equation}
  u(\vec x) = \exp(-10\norm{\vec x}^2).
\end{equation}
We discretize the problem given by \eqref{eq:Problem1} under the
algorithm set out in \S \ref{sec:discretisation}, numerical
convergence results are shown in Figure \ref{Fig:Nondiv:Prob1}.
\subsection{Test problem with convection dominated operator}
\label{test:advection-dominated}
The second test problem demonstrates the ability to overcome
oscillations introduced into the standard finite element when
rewriting the operator in divergence form. Take
\begin{gather}
  \label{eq:Problem2}
  a(\vec{x})
  = 
  \arctan{K(\norm{\vec x}^2-1)} + 2
  \\
  b(\vec x)
  = 0.
\end{gather}
with $K\in\reals^+$. Rewriting in divergence form gives
\begin{equation}
  \label{prob:divform}
    \frob{\A}{\Hess u} 
    =
    \div{\A \nabla u} 
    - 
    \div\A \nabla u.
\end{equation}
The derivatives
\begin{equation}
  \pd \alpha a(\vec x)
  = 
  \frac{d K x_\alpha} 
       {1 +K\left(\norm{\vec x}^2 - 1\right)}
\end{equation}
can be made arbitrarily large on the unit circle by choosing $K$
appropriately (see Figure ~\ref{Fig:atan}).

We choose our problem data $f$
such that the exact solution to the problem is given by:
\begin{equation}
  u(\vec x) = \sin{\pi x_1}\sin{\pi x_2}.
\end{equation}

We then construct the standard finite element method around
(\ref{prob:divform}), that is find $U\in \feszero$ such that for each
$\Phizero \in \feszero$
\begin{equation}
  \ltwop{\A\nabla U}{\nabla\Phizero} 
  -
  \ltwop{\div\A\nabla U}{\Phizero}
  =
  \ltwop{f}{\Phizero}.
\end{equation}
If $K$ is chosen small enough the standard finite element method
converges optimally. If we increase the value of $K$ oscillations
become apparent in the finite element solution along the unit
circle. Figure \ref{Fig:Fem:Prob2} demonstrates the oscillations
arising from this method compared to discretizing using the
nonvariational finite element method.

Figure \ref{Fig:Nondiv:Prob2} shows the numerical convergence rates of
the nonvariational finite element method applied to this problem.

\subsection{Test problem choosing a solution with nonsymmetric Hessian}
\label{test:nonsymmetric-hessian}
In this test we choose the operator such that $b(\vec x)$ is non-zero.
To maintain ellipticity in this problem we must choose $a(\vec x)$
such that the trace of $\vec A$ dominates it's determinant. We choose
\begin{gather}
  \label{eq:Problem3}
  a(\vec{x})
  = 
  2
  \\
  b(\vec x)
  = (x_1^2 x_2^2)^{1/3}.
\end{gather}
We choose the problem data such that the exact solution is given by
\begin{equation}
  u(\vec x) = 
  \begin{cases}
    \frac{x_1x_2(x_1^2-x_2^2)}{x_1^2+x_2^2}  &\geovec x \neq \geovec 0\\
    0 &\geovec x = \geovec 0.
  \end{cases}
\end{equation}
This function has a nonsymmetric Hessian at the point $\geovec 0$. The
nontrivial Dirichlet boundary is dealt with using Remark
\ref{rem:non-trivial-dirichlet}. Figure \ref{Fig:Nondiv:Prob3} shows
numerical results for this problem.

\subsection{Test problem with quasilinear PDE in nondivergence form}
\label{test:quasilin}
The problem under consideration in this test is the following
quasi-linear PDE arising from differential geometry:
\begin{equation}
  \label{eq:mcf}
  \div{\frac{\nabla u}{\sqrt{1+\norm{\nabla u}^2}}} 
  = 
  \frac{f}{\sqrt{1+\norm{\nabla u}^2}},
\end{equation}
where $\sqrt{1+\norm{\nabla u}^2}$ is the area element. Here we are
using $\norm{\nabla u}^2 = \D u \nabla u$. 
Applying a fixed point linearization given an initial guess $u^0$ for
each $n \in \naturals$ we seek $u^n$ such that
\begin{equation}
  \label{eq:mcflinear}
  \div{\frac{\nabla u^n}{\sqrt{1+\norm{\nabla u^{n-1}}^2}}} 
  = 
  \frac{f}{\sqrt{1+\norm{\nabla u^{n-1}}^2}}.
\end{equation}
Applying a standard finite element discretization of
(\ref{eq:mcflinear}) yields: Given $U^0\in\feszero$, for each $n \in
\naturals$ find $U^n\in\feszero$ such that for each $\Phizero\in\feszero$
\begin{equation}
  \ltwop{{\frac{\nabla U^n}{\sqrt{1+\norm{\nabla U^{n-1}}^2}}}}{\nabla\Phizero}
    = 
  \ltwop{\frac{f}{\sqrt{1+\norm{\nabla U^{n-1}}^2}}}{\Phizero}.
\end{equation}

In fact we can work on this problem combining the two nonlinear
terms. To do so we must first rewrite (\ref{eq:mcf}) into the form
$\frob{A(u,\nabla u)}{\Hess u} = f$.
\begin{equation}
  \begin{split}
    f 
    =& 
    \sqrt{1+\norm{\nabla u}^2}
    \div{\frac{\nabla u}
      {\sqrt{1+\norm{\nabla u}^2}}}
    \\
    =& 
    \sqrt{1+\norm{\nabla u}^2}
    \left(
    \frac{\Delta u}
         {\sqrt{1+\norm{\nabla u}^2}}
    +
    \frac{\D \left( 1 + \norm{\nabla u}^2 \right)}
         {2\left(1+\norm{\nabla u}^2\right)^{3/2}}
         \nabla u
         \right)
    \\
    =&
    \Delta u 
    +
    \frac{\D u \Hess u \nabla u}
         {1 + \norm{\nabla u}^2}
    \\
    =&
    \frob{\left( 
    \geovec I + \frac{\nabla u \D u}
        {1 + \norm{\nabla u}^2}
    \right)}
    {\Hess u}.       
  \end{split}
\end{equation}
Applying a similar fixed point linearization given an initial guess
$u^0$ for each $n \in\naturals$ we seek $u^n$ such that
\begin{equation}
  \label{eq:ndquasilinear}
  \frob{\left( 
    \geovec I + \frac{\nabla u^{n-1} \D u^{n-1}}
    {1 + \norm{\nabla u^{n-1}}^2}
    \right)}
       {\Hess u^n} 
       = f
\end{equation}

Discretizing the problem is then similar to that set out in Section
\ref{sec:discretisation}. The component matrixes $\nummat M$ and
$\nummat C_{\alpha,\beta}$ are problem independent, $\nummat
B^{\alpha,\beta}$ are defined as
\begin{equation}
  \nummat B^{\alpha,\beta}
  =
  \begin{cases}
    \ltwop
        {
          \numvec \Phizero
        }
        {
          1 + \frac
              {\pd \alpha U^{n-1} \pd \beta U^{n-1}}
              {1 + \norm{\nabla U^{n-1}}^2}
              \numvec \Phi}, 
        &\text{ for } \alpha = \beta,
        \\\\
    \ltwop
        {
          \numvec \Phizero
        }
        {
          \frac
              {\pd \alpha U^{n-1} \pd \beta U^{n-1}}
              {1 + \norm{\nabla U^{n-1}}^2}
              \numvec \Phi}, 
            &\text{ for } \alpha \neq \beta.
  \end{cases}
\end{equation}

Table \ref{tbl:quasilin} compares the two linearizations
(\ref{eq:mcflinear}) and (\ref{eq:ndquasilinear}). Figure
\ref{fig:quasilin} show asymptotic numerical convergence results
for NVFEM applied to (\ref{eq:ndquasilinear}).

\renewcommand{\figscale}{0.29}

\begin{figure}
  \caption{A visualization of the coefficient of the operators
    (\ref{eq:Problem1}) (on the left) and (\ref{eq:Problem2}) (on the
    right).
  }
  \begin{center}
    \subfigure[][{The function $(x_1^2 x_2^2)^{1/3} + 1$ over
        $\W$. Note the derivatives are singular at $x_1=0$ and
        $x_2=0$.}]{
      \label{Fig:x2y2}
      \includegraphics[scale=\figscale]{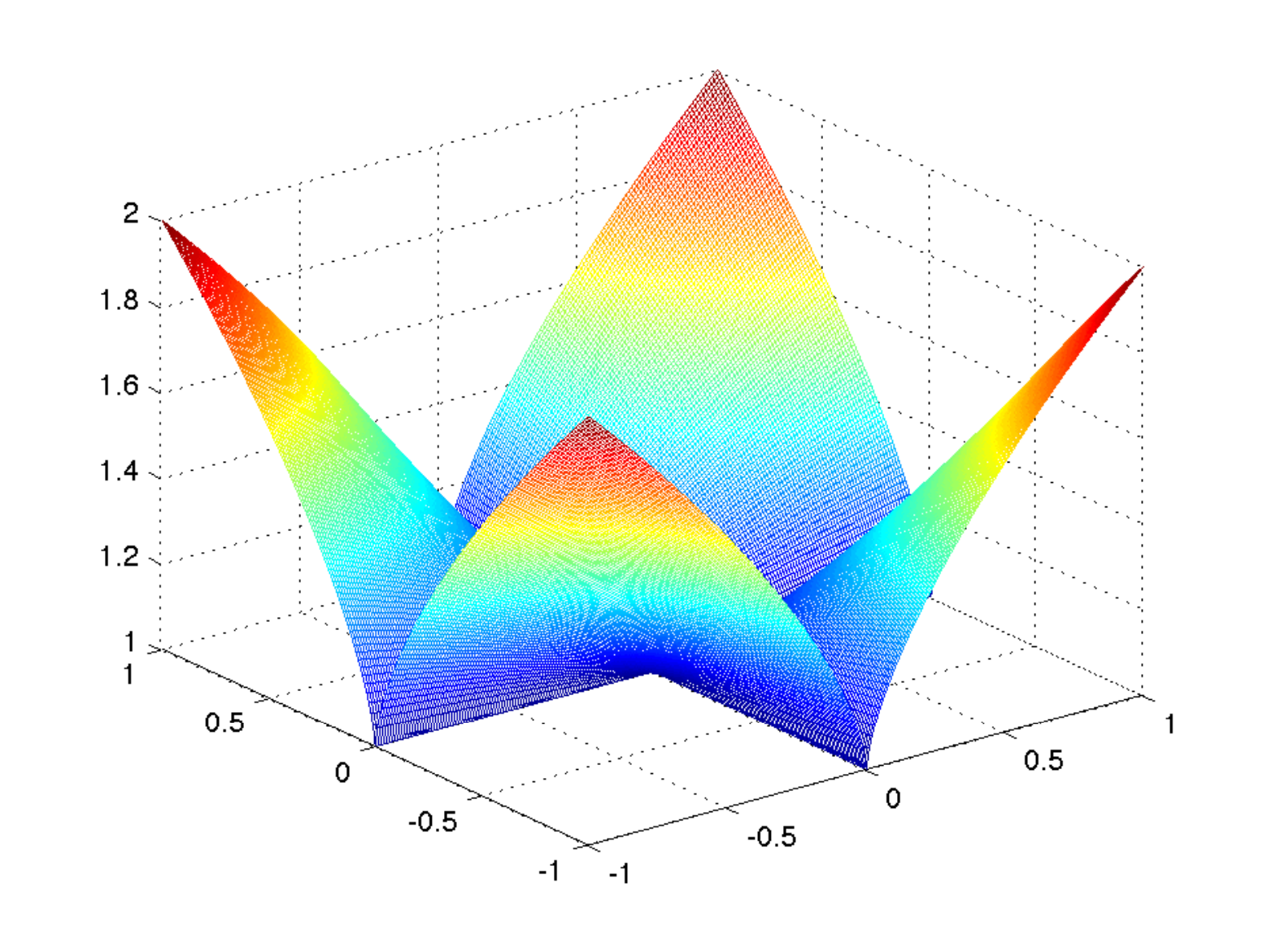}
    }
    \hfill
    \subfigure[][{The function $\arctan{5000(\norm{x}^2-1)}$ over
        $\W$. Note the derivatives are very large on the unit
        circle.}]{
      \label{Fig:atan}
      \includegraphics[scale=\figscale]{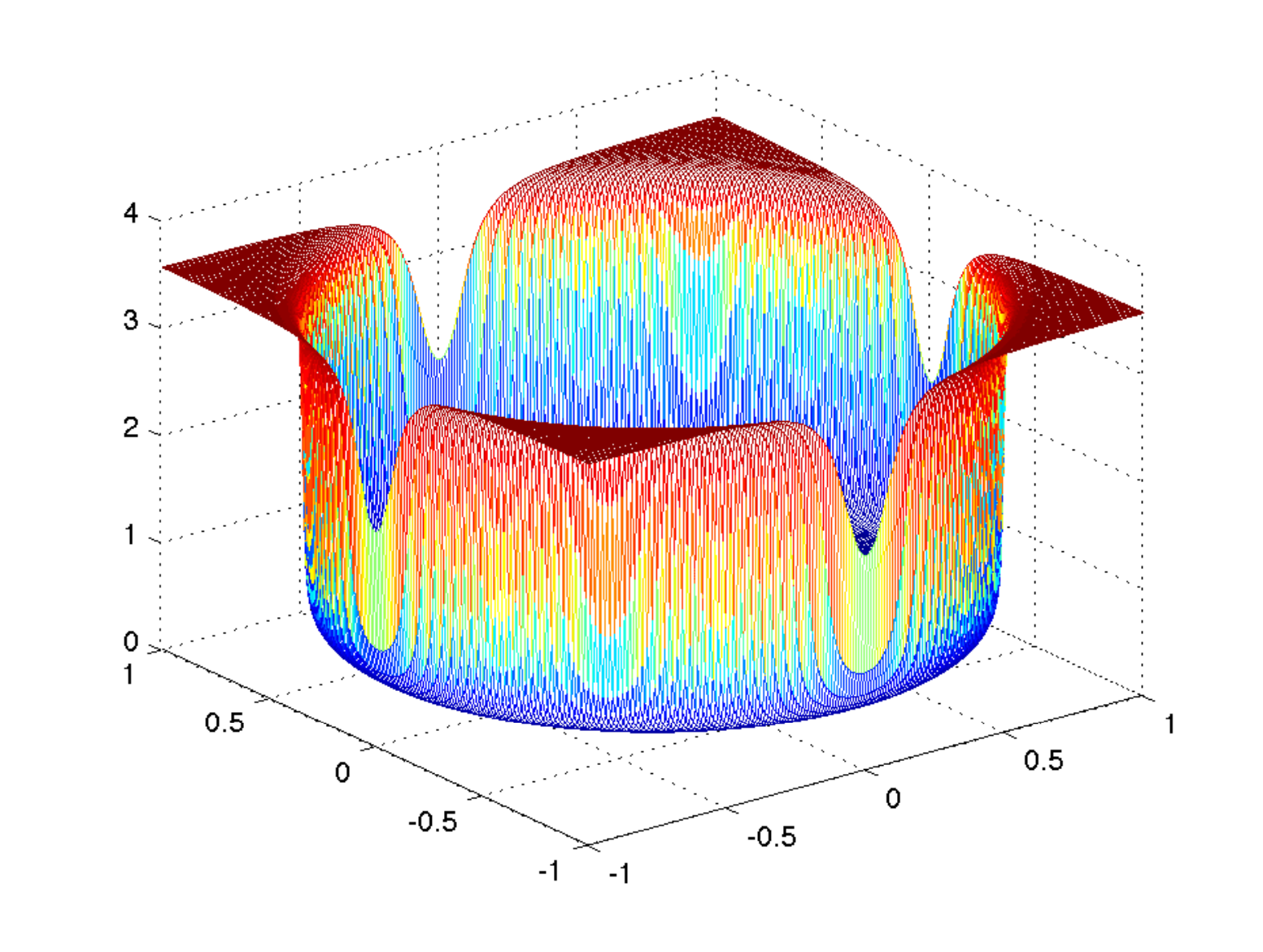}
    }
  \end{center}
\end{figure}

\renewcommand{\figscale}{0.4}

\begin{figure}[ht]
  \caption{  \label{Fig:Nondiv:Prob1}
    {Test \ref{test:nondiff}. Errors and convergence rates for the
      NVFEM applied to a non-divergence form operator
      \eqref{eq:Problem1}, choosing $f$ appropriately such that
      $u(\vec x) = \exp{(-10\norm{\vec x})}$. The convergence rates
      are optimal, that is for $\poly{1}$-elements (on the left)
      $\Norm{u - U} = \Oh(h^2)$ and $\norm{u - U}_1 = \Oh(h)$. For
      $\poly{2}$-elements (on the right) $\Norm{u - U} = \Oh(h^3)$ and
      $\norm{u-U}_1 = \Oh(h^2)$.}
    }
  \begin{center}
    \subfigure[][$\poly{1}$-elements]{
      \includegraphics[scale=\figscale]{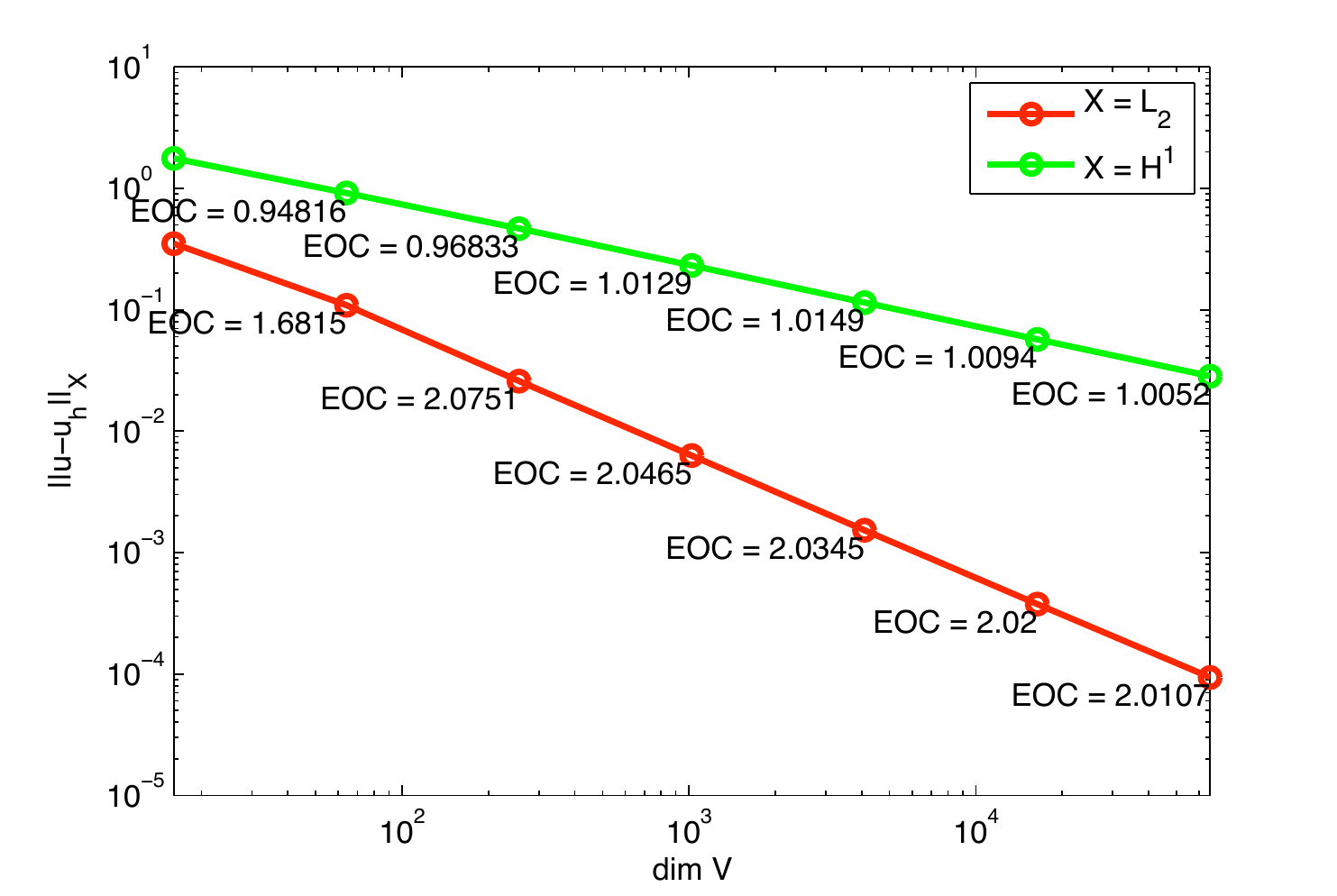}
    }
    \hfill
    \subfigure[][$\poly{2}$-elements]{
      \includegraphics[scale=\figscale]{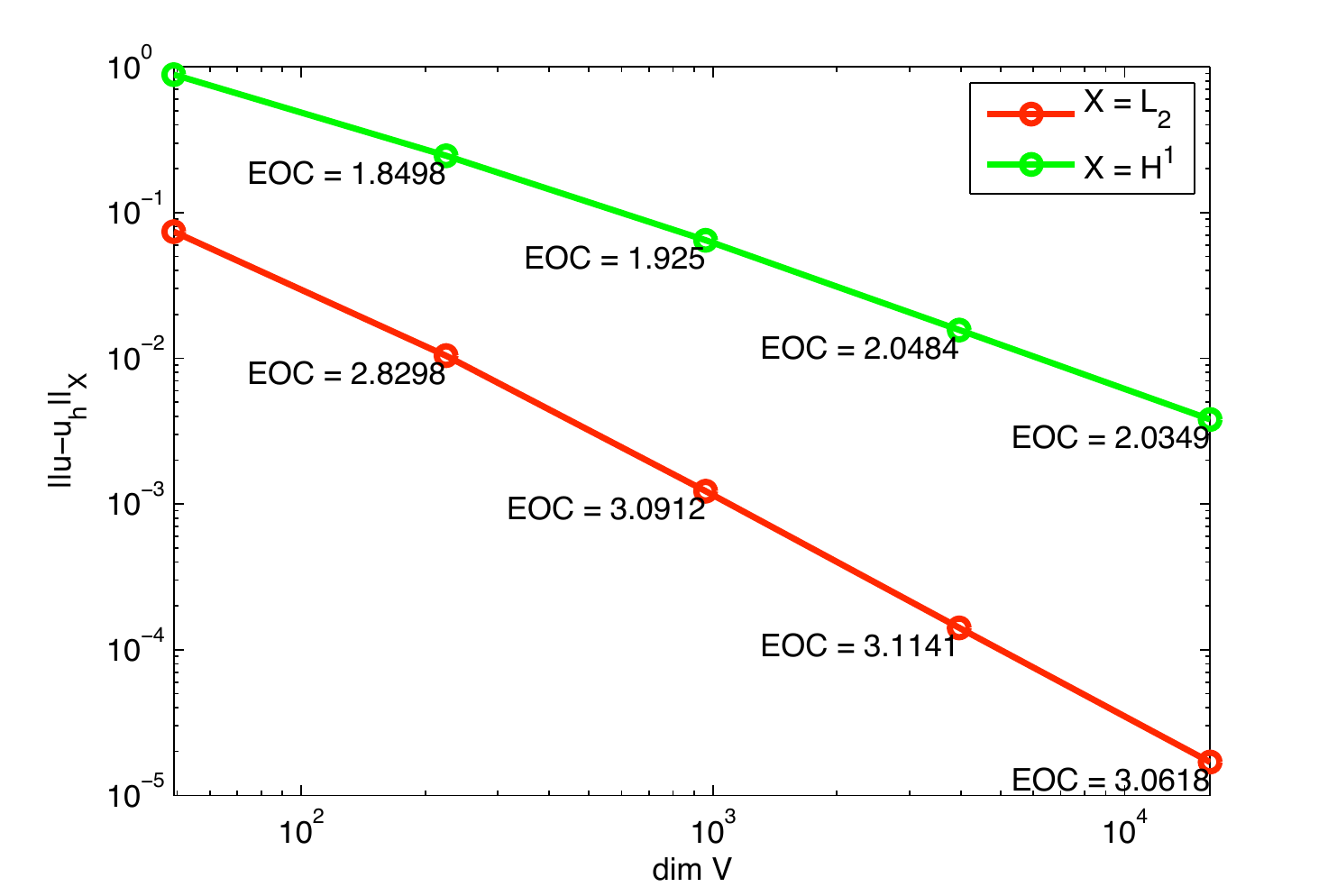}
    }
  \end{center}
\end{figure}

\begin{table}
  \caption{  \label{tbl:conditionnumber}
    Test \ref{test:nondiff}. On the condition number of $\nummat E$ upon
    discretizing problem \eqref{eq:Problem1} using $\poly{1}$ finite
    elements. As claimed in Remark \ref{rem:conditionnumber}
    $\kappa(\nummat E) \approx Ch^{-2}$. }
  \small
  \begin{tabular}{|c|c|c|c|}
    \hline 
    $\dim{\fes}$ & $h$ & $\kappa(\nummat E)$ & $h^{-2}\kappa(\nummat E)$\\ \hline 
    $16$ & $0.4714$ & $4.904 \times 10^1$ & $10.898$\\
    $64$ & $0.202$ & $6.594 \times 10^2$ & $26.952$\\
    $256$ & $0.0943$ & $3.665 \times 10^3$ & $32.633$\\  
    $1024$ & $0.0456$ & $1.722 \times 10^4$ & $35.833$\\ 
    $4096$ & $0.0224$ & $6.894 \times 10^4$ & $34.737$\\
    $16384$ & $0.0111$ & $3.383 \times 10^5$ & $41.949$\\
    $65536$ & $0.0055$ & $1.337 \times 10^6$ & $40.43$\\
    \hline
  \end{tabular}
  
\end{table}

\begin{figure}[ht]
  \caption{ \label{Fig:Nondiv:Prob2}
    {Test \ref{test:advection-dominated}. Errors and convergence rates
      for the NVFEM applied to a non-divergence form operator
      \eqref{eq:Problem2} with $K = 5000$, choosing $f$ appropriately such that
      $u(\vec x) = \sin{\pi x_1}\sin{\pi x_2}$. The convergence rates
      are optimal, that is for $\poly{1}$-elements (on the left)
      $\Norm{u - U} = \Oh(h^2)$ and $\norm{u - U}_1 = \Oh(h)$. For
      $\poly{2}$-elements (on the right) $\Norm{u - U} = \Oh(h^3)$ and
      $\norm{u-U}_1 = \Oh(h^2)$.}
    }

  \begin{center}
    \subfigure[][$\poly{1}$-elements]{
      \includegraphics[scale=\figscale]{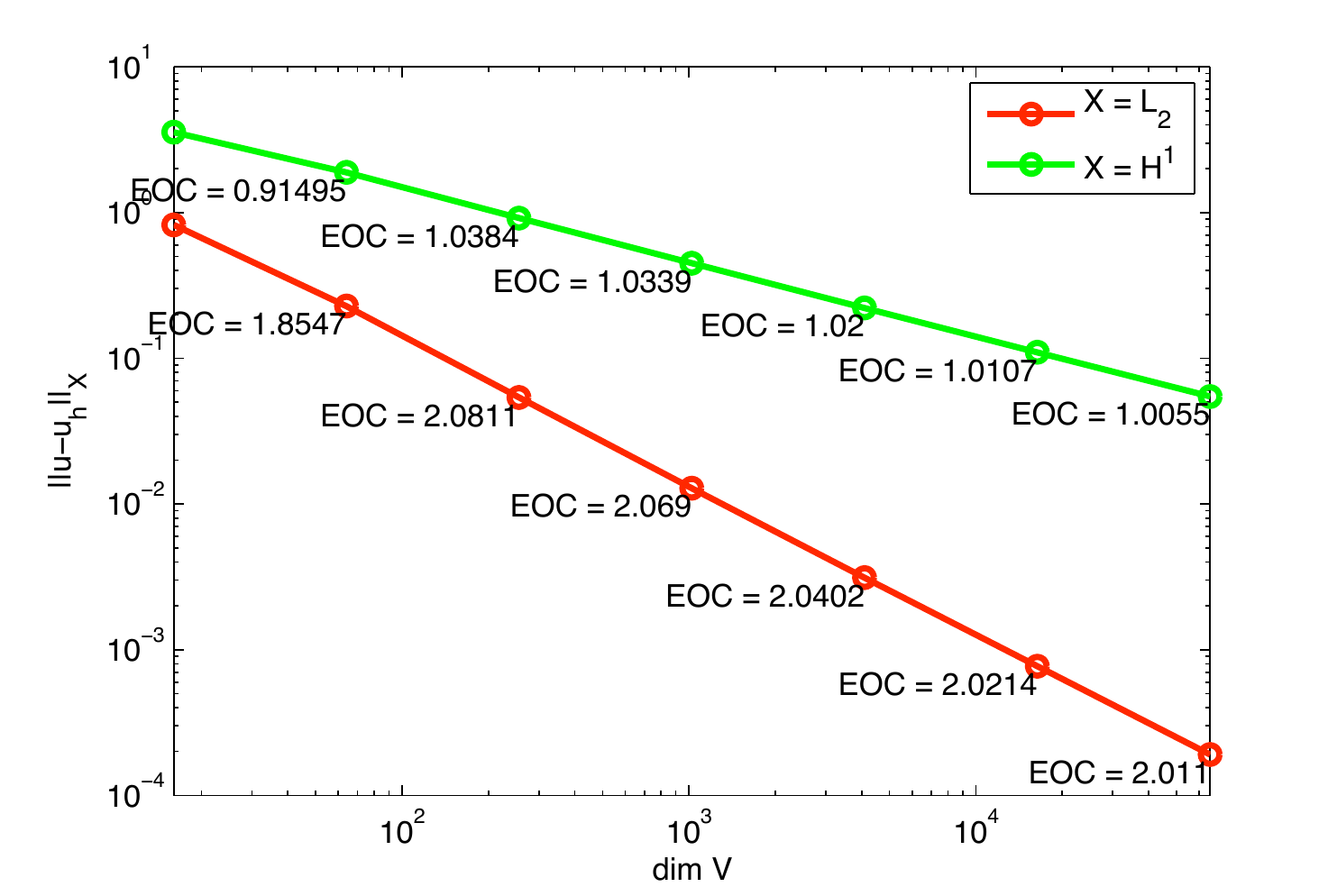}
    }
    \hfill
    \subfigure[][$\poly{2}$-elements]{
      \includegraphics[scale=\figscale]{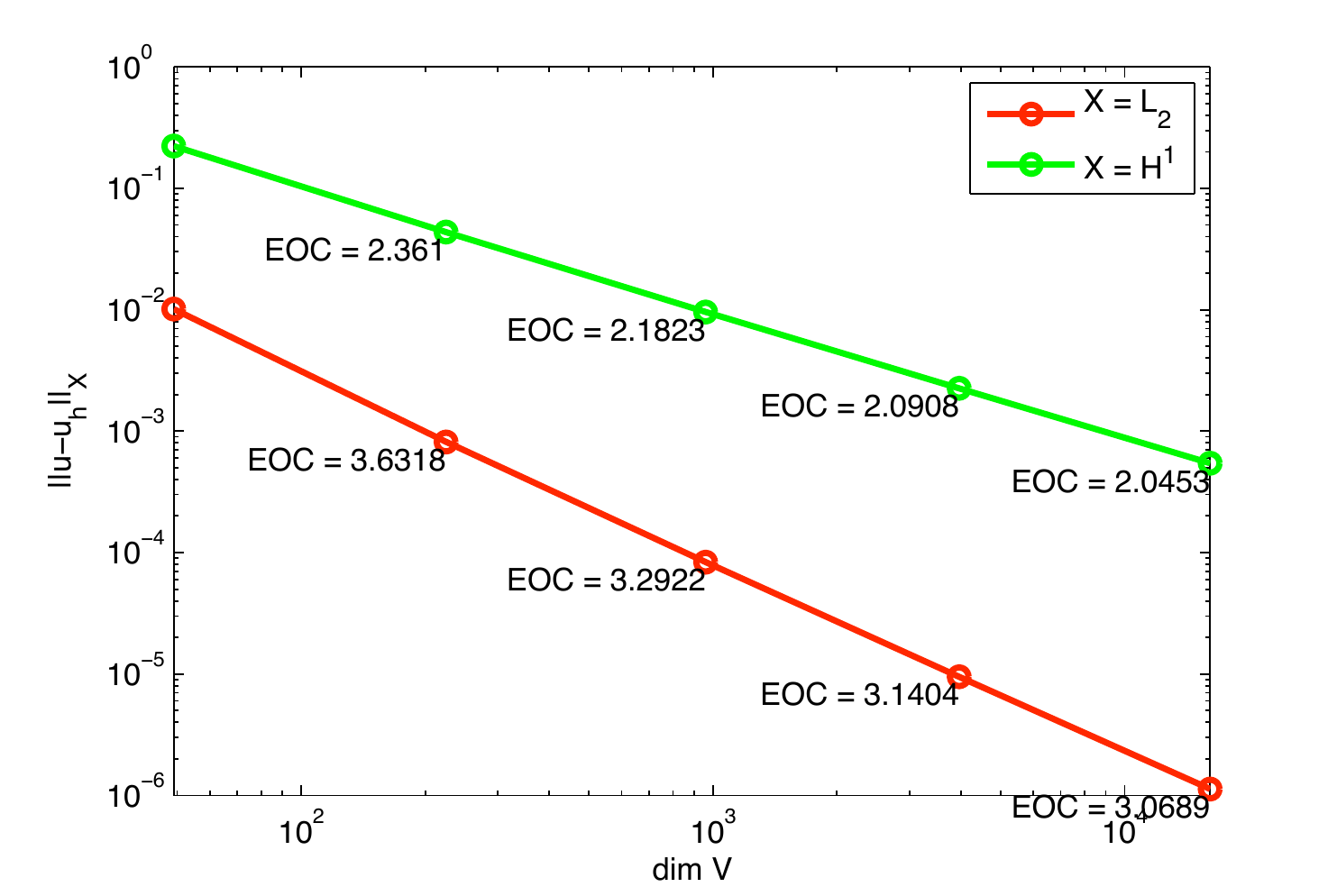}
    }
  \end{center}
\end{figure}

\renewcommand{\figscale}{0.65}

\begin{figure}[ht]
  \caption{ \label{Fig:Fem:Prob2} {Test
      \ref{test:advection-dominated}. On the left we present
      $\smash{\Norm{u - \tilde{U}}}_{\leb\infty(K)}$ plotted on a
      logarithmic scale as a function over $\W$. This represents the
      maximum error of the standard FE-solution, $\tilde{U}$, to
      problem (\ref{eq:Problem2})} with $16384$ DOF's ($h =
    1/32$). Notice the oscillations apparent on the unit circle.  On
    the right we show $\Norm{u - {U}}_{\leb\infty(K)}$ plotted on a
    logarithmic scale as a function over $\W$, the maximum error of
    the NVFE-solution, $U$, to problem (\ref{eq:Problem2}) with
    $16384$ DOF's ($h = 1/32$).}
  \begin{center}
    \includegraphics[scale=\figscale]{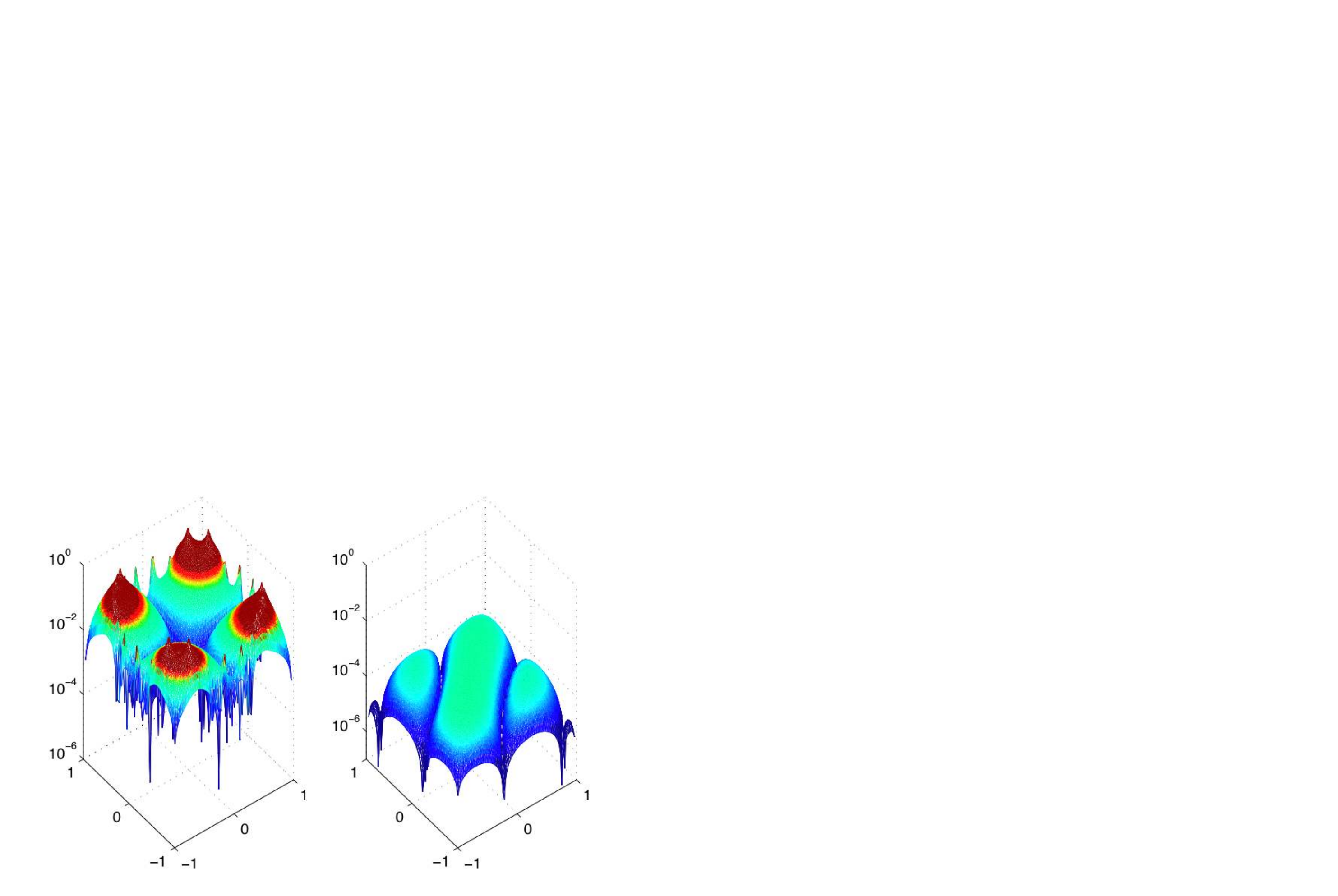}
  \end{center}
\end{figure}

\renewcommand{\figscale}{0.4}

\begin{figure}[ht]
  \caption{  \label{Fig:Nondiv:Prob3}
    {Test \ref{test:nonsymmetric-hessian}. Errors and convergence
      rates for the NVFEM on an operator (\ref{eq:Problem3}), choosing
      $f$ appropriately such that $u(\vec x) =
      \frac{x_1x_2(x_1^2-x_2^2)}{x_1^2+x_2^2}$ if $\vec x \neq \vec
      0$, or $u(\vec x) = 0$ otherwise. The convergence rates are
      optimal, that is for $\poly{1}$-elements (on the left) $\Norm{u
        - U} = \Oh(h^2)$ and $\norm{u - U}_1 = \Oh(h)$. For
      $\poly{2}$-elements (on the right) $\Norm{u - U} = \Oh(h^3)$ and
      $\norm{u-U}_1 = \Oh(h^2)$.}
  }
    \begin{center}
      \subfigure[][$\poly{1}$-elements]{
      \includegraphics[scale=\figscale]{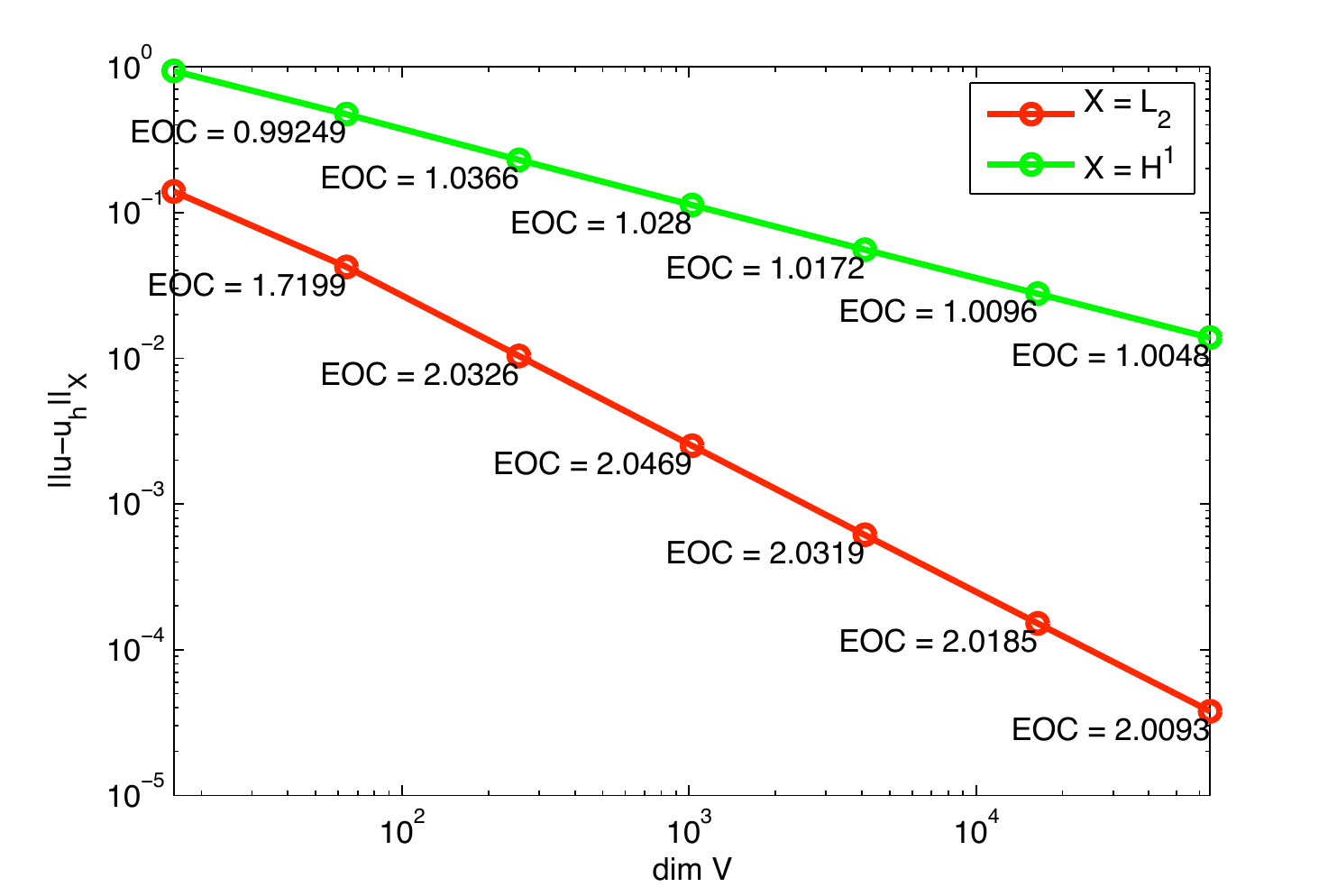}
      }
      \hfill
      \subfigure[][$\poly{2}$-elements]{
        \includegraphics[scale=\figscale]{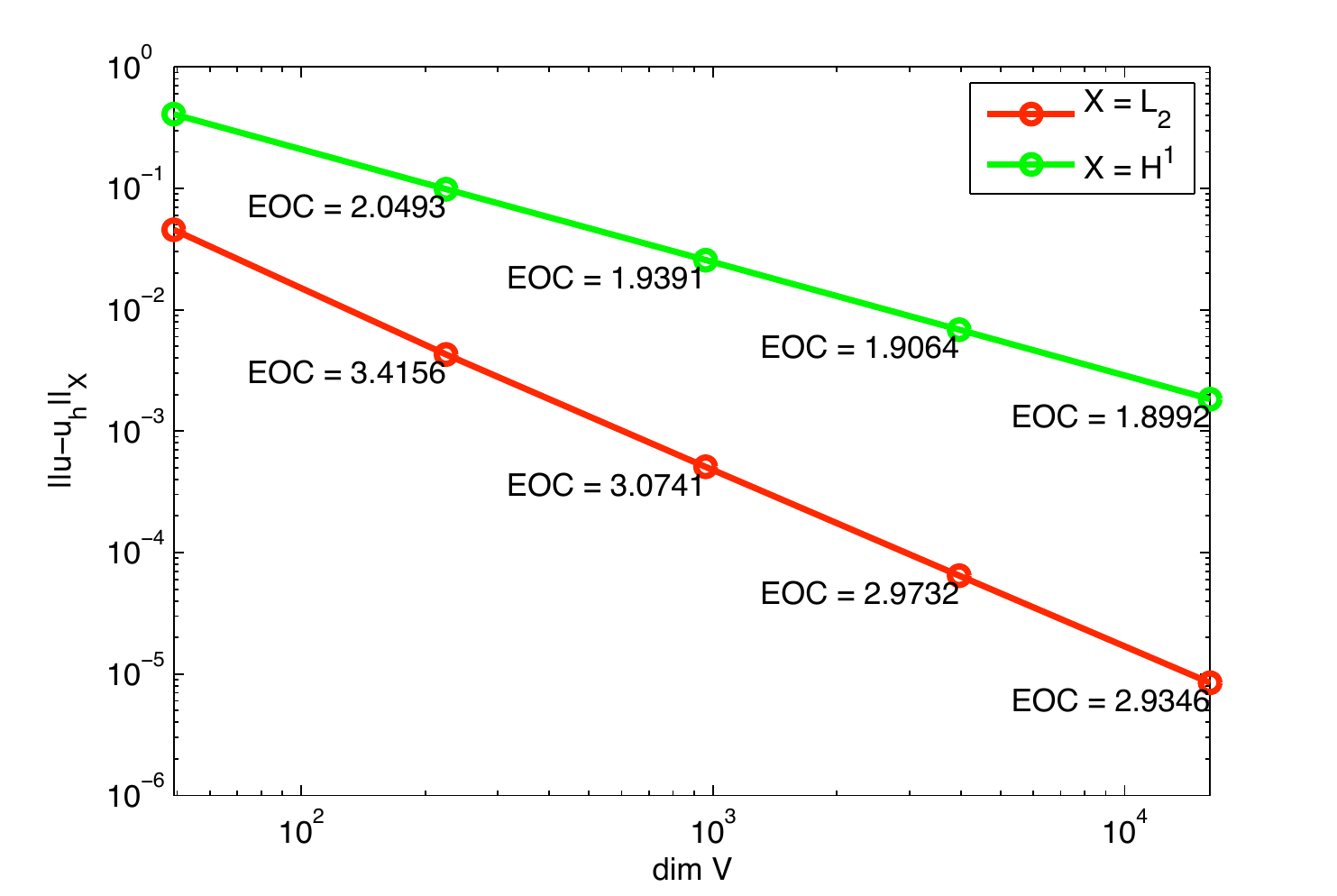}
        }
    \end{center}
\end{figure}

\renewcommand{\figscale}{0.29}

\begin{table}
  \caption{ \label{tbl:quasilin} Test \ref{test:quasilin}. Comparison
    of the fixed point linearization in variational form
    (\ref{eq:mcflinear}) and in nonvariational form
    (\ref{eq:ndquasilinear}). We fix $f$ appropriately such that
    $u(\vec x) = \sin{\pi x_1}\sin{\pi x_2}$. Taking initial
    guesses $U^0 = \tilde{U^0} = 0$ we discretize problem
    (\ref{eq:mcf}) using a standard FEM and using the NVFEM. Denoting
    $U_i$ and $\tilde{U}_i$ to be the NVFE-solution FE-solution
    respectively we run both linearizations for until a tolerance
    $\Norm{U_{n+1}-U_n}$
    (resp. $\smash{\Norm{\tilde{U}_{n+1}-\tilde{U}_n}}$) $\leq h^2$ is
    achieved. We compute both the stagnation point---which is the
    iteration at which the prescribed tolerance is achieved---and the
    total CPU time. Notice there is significant savings in the number
    of iterations required to reach the stagnation point using the
    NVFEM over the standard FEM, however each iteration is
    computationally more costly using the NVFEM since the system is
    larger and more complicated to solve. The CPU cost for the entire
    algorithm is comparable for each fixed $h$.}
\small
\begin{tabular}{|c|c|c|c|c|c|c|c|}
  \hline 
  & $h$  & $\sqrt{2}/5$ & $\sqrt{2}/10$ & $\sqrt{2}/20$ & $\sqrt{2}/40$ & $\sqrt{2}/80$ & $\sqrt{2}/160$\\ 
  \hline 
      {FEM} 
      & Stag. Point& $5$ & $13$ & $16$ & $26$ & $32$ & $36$ 
      \\
      & CPU Time& $0.50$ & $4.02$ & $17.51$ & $117.58$ & $796.58$ & $5308.81$
      \\
      \hline
          {NDFEM}
          & Stag. Point& $4$ & $6$ & $7$ & $8$ & $10$ & $12$
          \\
          & CPU Time& $0.72$ & $3.40$ & $16.49$ & $97.93$ & $838.8$ & $5256.84$
          \\
          \hline
\end{tabular}

\end{table}

\renewcommand{\figscale}{0.4}

\begin{figure}[ht]
  \caption{\label{fig:quasilin}
    {Test \ref{test:quasilin}. Errors and convergence rates for NVFEM
      applied to (\ref{eq:mcf}), a quasi-linear PDE under a fixed
      point linearization. We fix $f$ appropriately such that
      $u(\vec x) = \sin{\pi x_1}\sin{\pi x_2}$, taking an initial
      guess $u^0 = 0$. The convergence rates are optimal, that is for
      $\poly{1}$-elements (on the left) $\Norm{u - U} = \Oh(h^2)$ and
      $\norm{u - U}_1 = \Oh(h)$. For $\poly{2}$-elements (on the
      right) $\Norm{u - U} = \Oh(h^3)$ and $\norm{u-U}_1 =
      \Oh(h^2)$.}}
  \begin{center}
    \subfigure[][$\poly{1}$-elements]{
      \includegraphics[scale=\figscale]{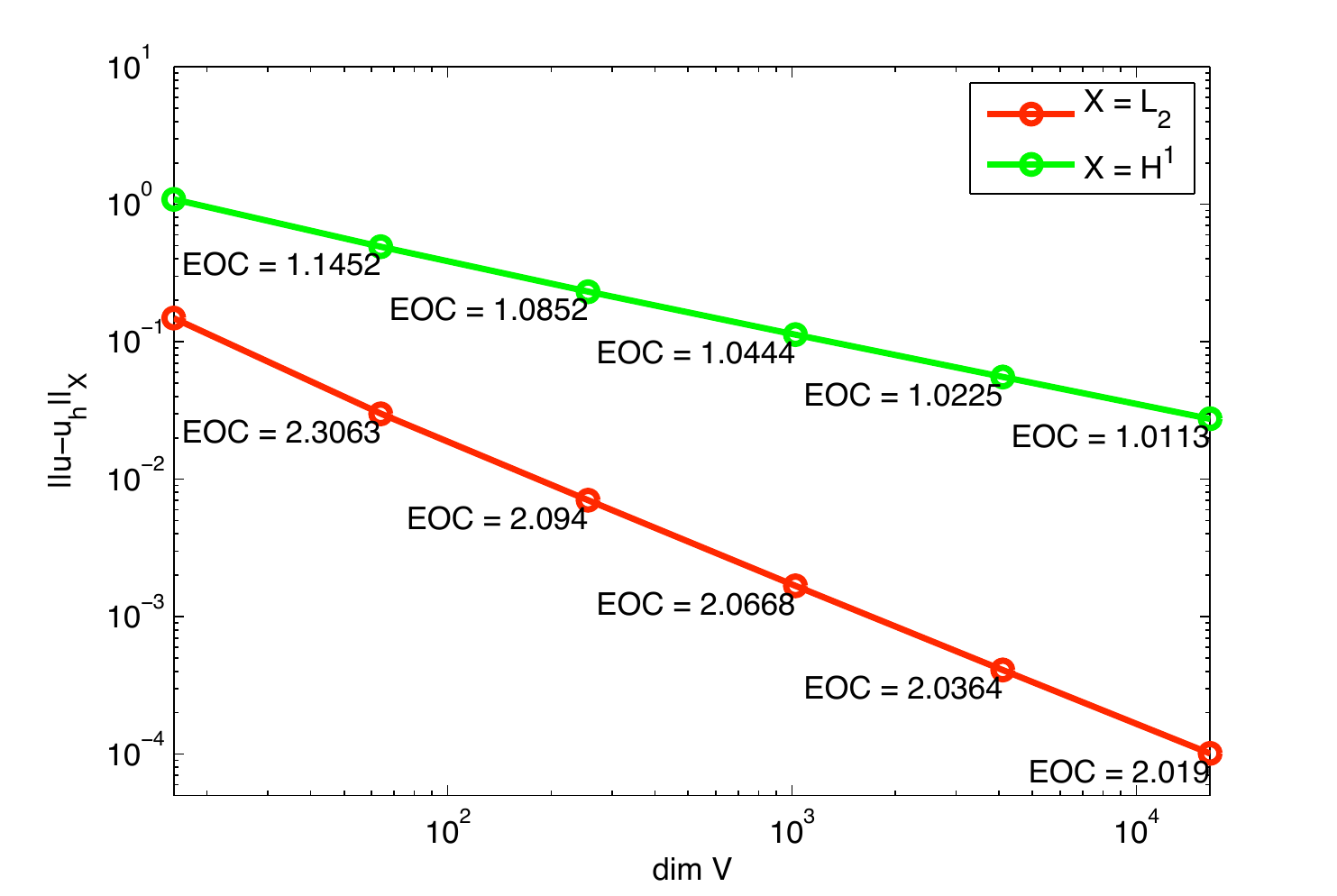}
    }
    \hfill
    \subfigure[][$\poly{2}$-elements]{
      \includegraphics[scale=\figscale]{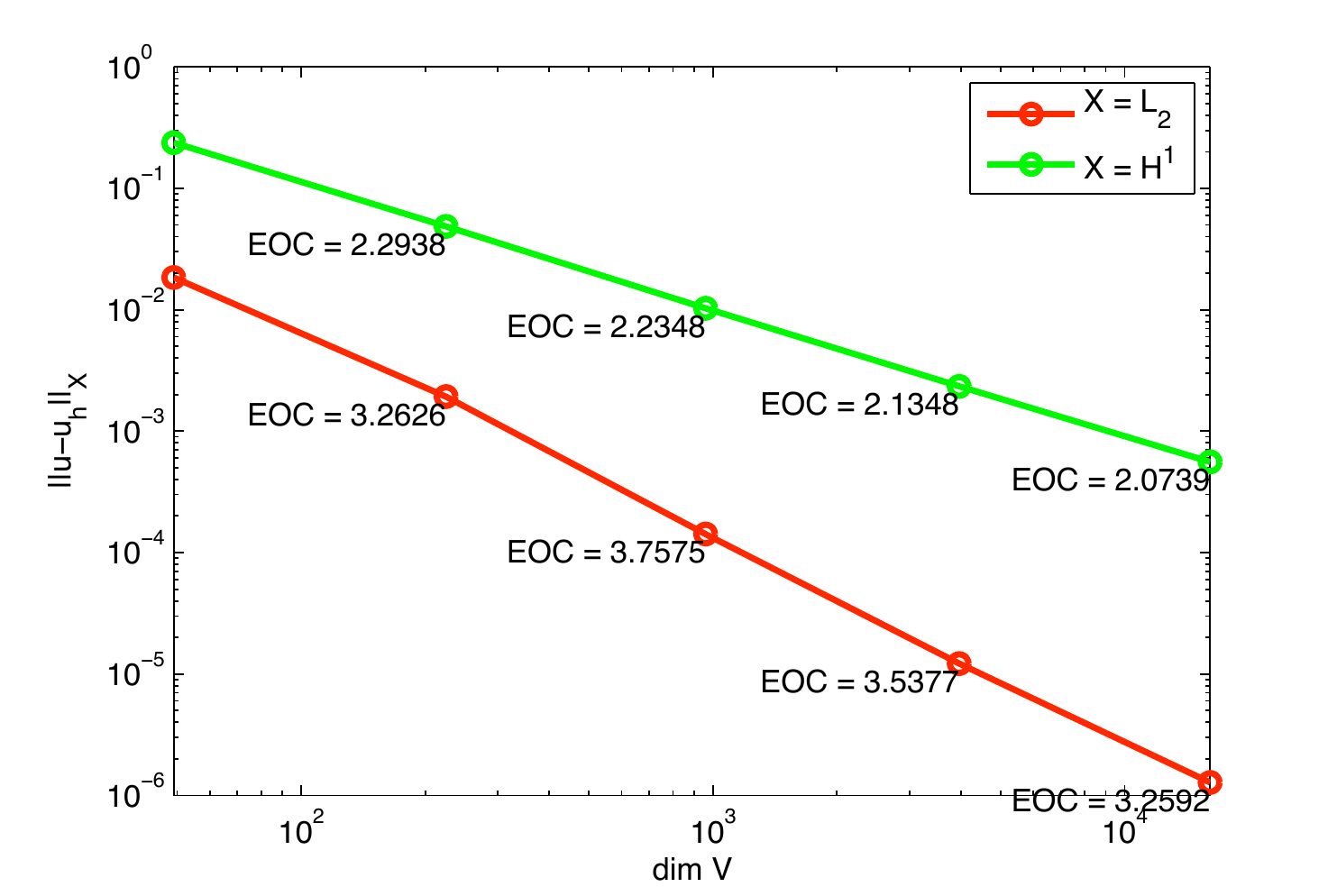}
    }
  \end{center}
\end{figure}



\newcommand{\etalchar}[1]{$^{#1}$}


\end{document}